# Approximation by normal distribution for a sample sum in sampling without replacement from a finite population

Ibrahim Bin Mohamed,

University of Malaya, Malaysia,

imohamed@um.edu.my

Sherzod M. Mirakhmedov,

Institute of Mathematics, Tashkent.

shmirakhmedov@yahoo.com

**Abstract**.

A sum of observations derived by a simple random sampling design from a population of independent random variables is studied. A procedure finding a general term of Edgeworth asymptotic expansion is presented. The Lindeberg condition of asymptotic normality, Berry-Esseen bound , Edgeworth asymptotic expansions under weakened conditions and Cramer type large deviation results are derived.

**1. Introduction**.

Let $\pi_N = (Y_{1N},...,Y_{NN})$ be a population of independent random variables (r.v.) and $r = (r_1,...,r_N)$ be a random vector (r.vec.) independent of $Y_{1N},...,Y_{NN}$ and such that $P\{r_1 = k_1,...,r_N = k_N\} = (N!)^{-1}$, for any permutation $k = (k_1,...,k_N)$ of the numbers $1,...,N$. Let $S_{nN} = Y_{r_1 N} + ... + Y_{r_n N}$, that is a sum of $n$ r.v.s chosen at random without replacement from the population $\pi_N$. We are interested in the approximation of the distribution of $S_{n,N}$ by a normal distribution.

The sample sum $S_{nN}$ is of interest in the context of its statistical applications. Under the hypothesis of homogeneity of two samples, the two-sample linear rank statistics can be reduced to that of sample sum $S_{nN}$ when the elements of the population are non-random real numbers; in this case, the elements of the population become the scores of the rank statistic, see Bickel and von Zwet (1978). Situation where $Y_{kN}$ are random variables arises, for instance, in the problems associated with two-stage samples, in particular in stratified random sampling, in a linear estimation for the mean total of a stratified population, see Cohran (1963). Note also if $n = N$, then $S_{nN}$ simply is a sum of independent r.v.s. Many authors have studied the sample sum. In the situation when $Y_{k,N}$ are non-random a sufficient and necessary condition, similar to the well-known Lindeberg condition, for asymptotic normality of $S_{nN}$ have been presented by Erdös and Réniy (1959) and Hajek (1960); we should refer also to Wald and Wolfowitz (1944) who for the first time have proved asymptotic normality of $S_{n,N}$ under certain conditions. The rate in the central limit theorem (CLT) was obtained by Bikelis (1969) and Höglund (1978). The second order approximation results was derived by Robinson (1978), Bickel and van Zwet (1978), Babu and Singh (1985), Babu and Bai (1996) and Bloznelis (2000). Also Cramer's large deviation result has been proved by Robinson (1977) and Hu et al (2007b). The case when the elements of the population $\pi_N$ are r.v.s has been extensively studied since the paper of von Bahr (1972), who have established for the first time a bound for the remainder term in CLT; later non-uniform variant of von Bahr's result was obtained by Mirakhmedov and Nabiev (1976). Von Bahr's bound has been improved in Corollary of Mirakhmedov (1985) and in Theorem 1 of Zhao et al (2004). Two terms Edgeworth asymptotic expansion result has been established by Mirakhmedov (1979) and Hu et al (2007a), whereas expansion with arbitrary number of terms was obtained by Mirakhmedov (1983). Large deviation results follows from Theorems 11 and 12 of Mirakhmedov (1996).

In this paper we consider the general situation when the elements of the population $\pi_N$ are r.v.s which may be degenerate (i.e. non-random). The aim of the paper is threefold. First, in Sec 2 we present a procedure finding a general term of Edgeworth asymptotic expansion for the sample sum; our rule based on slightly generalized formula of Erdos and Reny (1959), and used fact that the integrant of that formula is the characteristic function (ch.f.) of an independent two-dimensional r.vec.s. This rule is considerably simpler than that of Mirakhmedov (1983), who have used formula of von Bahr (1972), and Bloznelis (2000), who gives just two terms of asymptotic expansion for the case when the elements of a population are non-random numbers. Second, we obtain Berry-Esseen bound and Edgeworth expansion result under weakened moment conditions, namely we assume that $E|Y_{mN}|^{k+1+\delta} < \infty$, $m = 1, ..., N$, where $0 < \delta \leq 1$, $k$ is the number of terms in asymptotic expansion (notice that the normal distribution is the first term). A general Theorems 3.1 and 3.2 of Sec 3 gives the bound for the remainder term in a unusually form; this together with above mentioned fact on the integrand of Erdos and Reny formula allows to write estimations of the remainder terms in term of Lyapunov ratios of a associated r.v.s $Z_m$, see (2.2) and (2.5). By-turn this fact allows us to improve concerning theorems of Zhao et al (2004) and Hu et al (2007a) by weakening the moment conditions, bringing some correction and giving one more term of expansion, see Remarks 3.4 and 3.5. Third, we note that the sample sum is a special case of so namely decomposable statistics, considered for instance in Mirakhmedov et al (2012). Using this fact we derive for the first time the Lindeberg type condition of asymptotic normality and Cramer's type large deviation result. These results are extension of concerning theorems of Hajek (1969) and Robinson (1977) respectively. Thus in the paper complete spectrum of results on approximation by normal distribution is collected. Main assertions are presented in Sec 3, whereas the proofs are given in Section 4.

In what follows $Bi(p)$ is used to denote Bernoulli distribution with probability of success $p$; $\mathcal{L}(\xi)$ stands for the distribution of the r.v. (r.vec.) $\xi$; $c$, $C$ with or without an index are positive universal constants whose value may differ at each occurrence. In the sequel notations we shall suppress the dependence on $N$ whenever it is convenient.

**2. Procedure of finding a general term of asymptotical expansion of the ch.f. of $S_{nN}$.**

Write $S_{nN} = Y_{1N}\eta_1 + ... + Y_{NN}\eta_N$, where $\eta_m = 1$ if $Y_{mN}$ appears in the sample else $\eta_m = 0$, i.e. $\mathcal{L}(\eta_m) = Bi(p)$ with $p = n/N$ and $\eta_1 + ... + \eta_N = n$. The joint distribution of $(\eta_1, ..., \eta_N)$ can be viewed as conditional joint distribution of the r.vec. $(\xi_1, ..., \xi_N)$ under $\zeta_N := \xi_1 + ... + \xi_N = n$, where $\xi_1, ..., \xi_N$ are independent r.v.s and $\mathcal{L}(\xi_m) = Bi(p)$. Hence $E\left(e^{it\bar{S}_{nN}} / \zeta_N = n\right) = Ee^{itS_{nN}}$, where $\bar{S}_{nN} = Y_{1N}\xi_1 + ... + Y_{NN}\xi_N$. This together with

$$E\left(e^{it\bar{S}_{nN} + i\tau(\zeta_N - n)}\right) = E\left(e^{i\tau(\zeta_N - n)} E\left(e^{it\bar{S}_{nN}} / \zeta_N\right)\right) = \sum_{k=0}^{\infty} e^{i\tau(k-n)} P\{\zeta_N = k\} E\left(e^{it\bar{S}_{nN}} / \zeta_N = k\right).$$

and fact that $\hat{S}_{nN}$ and $\zeta_N$ are a sums of independent r.v.s, implies by Fourier inversion

$$E\exp\{itS_{nN}\} = \frac{1}{\sqrt{2\pi}\,d_n(p)} \int_{-\pi\sqrt{nq}}^{\pi\sqrt{nq}} \prod_{m=1}^{N} E\exp\{itY_{mN}\xi_m + i\tau(\xi_m - p)/\sqrt{nq}\}d\tau, \qquad (2.1)$$

where $q = 1-p$, $d_n(p) := \sqrt{2\pi nq}\,P\{\zeta_N = n\} = \sqrt{2\pi nq}\,C_N^n p^n q^{N-n}$.

If $Y_{mN} = a_{mN}$, $m = 1,...,N$ are non-random real numbers then (2.1) can be written in the form

$$E\exp\{itS_{nN}\} = \frac{1}{\sqrt{2\pi}\,d_n(p)} \int_{-\pi\sqrt{nq}}^{\pi\sqrt{nq}} e^{-i\tau\sqrt{n/q}} \prod_{m=1}^{N}\left(q + p\exp\{ita_{mN} + i\tau/\sqrt{nq}\}\right)d\tau.$$

In fact this is the formula of Erdös and Réniy (1959) which had been used to prove the first satisfactory CLT for the sample sum. The formula had motivated for a number of related studies; for instance Hajek (1960), Bikelis (1969), Babu and Singh (1985) and Bloznelis (2000). Set

$$\gamma = N^{-1}\sum_{m=1}^{N} EY_{mN}, \quad \sigma^2 = \frac{1}{N}\sum_{m=1}^{N}\left[E(Y_{mN} - \gamma)^2 - p\left(E(Y_{mN} - \gamma)\right)^2\right],$$

$$Z_{mN}^x = (Y_{mN} - \gamma)x - pE(Y_{mN} - \gamma). \qquad (2.2)$$

Note that $ES_{nN} = n\gamma$ and $S_{nN} - n\gamma = Z_{1N}^{\eta_1} + ... + Z_{NN}^{\eta_N}$. From formula (2.1) we obtain

$$\varphi_n(t) := E\exp\left\{it\frac{S_{nN} - n\gamma}{\sigma\sqrt{n}}\right\} = \frac{\Theta_N(t)}{d_N(p)} = \frac{\Theta_N(t)}{\Theta_N(0)}, \qquad (2.3)$$

where

$$\Theta_N(t) = \frac{1}{\sqrt{2\pi}} \int_{-\pi\sqrt{nq}}^{\pi\sqrt{nq}} \prod_{m=1}^{N} \psi_m(t,\tau)d\tau, \qquad (2.4)$$

$$\psi_m(t,\tau) = E\exp\left\{\frac{itZ_{mN}^{\xi_m}}{\sigma\sqrt{n}} + \frac{i\tau(\xi_m - p)}{\sqrt{nq}}\right\},$$

because $d_n(p) = \Theta_N(0)$, due to inversion formula for the local probability. In what follows just to keep notation simple we put

$$Z_m = Z_{mN}^{\xi_m}. \qquad (2.5)$$

Note that $n\sigma^2 = VarZ_1 + ... + VarZ_N$, and

$$\sum_{m=1}^{N}\text{cov}(Z_m, \xi_m) = 0. \qquad (2.6)$$

In regard to formula (2.3) we remark that the integrand in $\Theta_N(t)$ is the ch.f. of a sum of independent two dimensional r.vec.s. This fact is crucial in our next considerations, in particular, in constructing the terms of asymptotic expansion of the ch.f. of the sample sum. Indeed, let $E|Z_m|^k < \infty$, $k \geq 3$, and $P_{m,N}(t,\tau)$, $m = 1,2,...$, be the well-known polynomials of $\tau$ and $t$ from the theory of asymptotic expansion of the ch.f. of a sum of independent r.vec.s , see Bhattacharya and

Rao (1976) p.52 (henceforth to be referred to as BR). In our case, the concerning sum is
$(Z,\xi) = (\tilde{Z}_1, \tilde{\xi}_1) + ... + (\tilde{Z}_N, \tilde{\xi}_N)$, where

$$\tilde{Z}_m = Z_m / \sigma\sqrt{n} \text{ and } \tilde{\xi}_m = (\xi_m - p)/\sqrt{nq}. \tag{2.7}$$

It is essential that this sum has *zero* expectation, *uncorrelated* components and a *unit* covariance matrix, due to (2.6). Hence the integrand in $\Theta_N(t)$, being the ch.f. of the sum $(Z,\xi)$, can be approximated by a power-series in $N^{-1/2}$ whose coefficients are polynomials of $t$ and $\tau$ containing the common factor $\exp\{-(t^2 + \tau^2)/2\}$, hence the series can be integrated wrt $\tau$ over the interval $(-\infty, \infty)$. The degree of $P_{m,N}(t,\tau)$ is $3m$ and minimal degree is $m+2$; the coefficients of $P_{m,N}(t,\tau)$ only involve the cumulants of the random vectors $(Z,\xi)$ of order less than or equal to $m+2$, in particular:

$$\frac{1}{\sqrt{N}} P_{1,N}(t,\tau) = \frac{i^3}{6} \sum_{m=1}^{N} E\left(t\tilde{Z}_m + \tau\tilde{\xi}_m\right)^3, \tag{2.8}$$

$$\frac{1}{N} P_{2,N}(t,\tau) = \frac{i^4}{24} \sum_{m=1}^{N} \left( E\left(t\tilde{Z}_m + \tau\tilde{\xi}_m\right)^4 - 3\left(E\left(t\tilde{Z}_m + \tau\tilde{\xi}_m\right)^2\right)^2 \right) + \frac{1}{2N} P_{1,N}^2(t,\tau). \tag{2.9}$$

Define functions $G_{m,N}(t)$ such that $G_{0,N}(t) = 1$ and

$$G_{m,N}(t) = \frac{1}{\sqrt{2\pi}} \int_{-\infty}^{\infty} P_{m,N}(t,\tau) \exp\left\{-\frac{\tau^2}{2}\right\} d\tau, \ m=1,2,...,$$

and write

$$Q_{k,N}(t) = e^{-\frac{t^2}{2}} \sum_{m=0}^{k} N^{-m/2} G_{m,N}(t). \tag{2.10}$$

Due to Theorem 9.11 of BR (1976) this series is $(k+1)$-term asymptotic expansion of the integral $\Theta_N(t)$. Finally, we get the $(k+1)$-term asymptotic expansion of $\varphi_n(t)$ by dividing $Q_{k,N}(t)$ with $Q_{k,N}(0)$. Also in order to estimate remainder terms in this approximation one can use directly Theorem 9.11, Lemma 9.5 and Lemma 14.3 of BR in an evidently simplified form. In particular, three term asymptotic expansion of the $\varphi_n(t)$ is

$$W_{3,N}(t) = e^{-\frac{t^2}{2}} \left(1 + N^{-1/2} G_{1,N}(t) + N^{-1}(G_{2,N}(t) - G_{2,N}(0))\right)$$

$$= e^{-t^2/2} \left\{ 1 + \frac{(it)^3}{6} \sum_{m=1}^{N} E\tilde{Z}_m^3 + \frac{(it)^6}{72} \left(\sum_{m=1}^{N} E\tilde{Z}_m^3\right)^2 \right.$$

$$\left. + \frac{(it)^4}{24} \left( \sum_{m=1}^{N} E\tilde{Z}_m^4 - 3\left( \sum_{m=1}^{N} \left(E\tilde{Z}_m^2\right)^2 + \left(\sum_{m=1}^{N} E\tilde{Z}_m^2 \tilde{\xi}_m\right)^2 \right) \right) + (it)^2 \frac{pq\alpha_{02}}{2n\sigma^2} \right\}$$

$$= e^{-t^2/2} \left\{ 1 + \frac{(it)^3}{6\sqrt{n}\sigma^3} \Lambda_1 + \frac{(it)^6}{72n\sigma^6} \Lambda_1^2 + \frac{(it)^4}{24n\sigma^4} \Lambda_2 + (it)^2 \frac{pq\alpha_{02}}{2n\sigma^2} \right\}, \tag{2.11}$$

where

$$\alpha_{kl} = \alpha_{kl}^{(1)}, \quad \alpha_{kl}^{(i)} = \frac{1}{N} \sum_{m=1}^{N} \left(E(Y_{mN} - \gamma)^k\right)^i \left(E(Y_{mN} - \gamma)\right)^l. \tag{2.12}$$

$$\Lambda_1 = \alpha_{30} - 3p\alpha_{21} + 2p^2\alpha_{03}, \quad \Lambda_2 = \alpha_{40} - 4p\alpha_{31} + 12p^2\alpha_{22} - 6p^3\alpha_{04} - 3p\alpha_{20}^{(2)} - 3q(\alpha_{20} - 2p\alpha_{02})^2.$$

The terms of asymptotic expansion of the distribution function of $S_{nN}/\sigma\sqrt{n}$ can be obtained by standard way: formally substitute

$$(-1)^v \frac{d^v}{du^v} \Phi(u) = -e^{-u^2/2} H_{v-1}(u)/\sqrt{2\pi}, \text{ where } \Phi(u) = \frac{1}{\sqrt{2\pi}} \int_{-\infty}^{u} e^{-\frac{t^2}{2}} dt,$$

instead of $(it)^v e^{-t^2/2}$ for each v in the expression for $W_{3,N}(t)e^{-t^2/2}$, see Lemma 7.2 of BR, p. 53, where $H_v(x)$ is the v-th order Hermit-Chebishev polynomial. Denoting $\mathbb{W}_{3,N}(u)$ for three term asymptotic expansion we have

$$\mathbb{W}_{3,N}(u) = \Phi(u) - \frac{e^{-u^2/2}}{\sqrt{2\pi}} \left\{ \frac{H_2(u)}{6\sigma^3 \sqrt{n}} \Lambda_1 \right.$$

$$\left. + \frac{1}{72n\sigma^6} \left[ H_5(u)\Lambda_1^2 + 3H_3(u)\sigma^2 \Lambda_2 + 36H_1(u)pq\sigma^4 \alpha_{02} \right] \right\}, \tag{2.13}$$

with $H_1(u) = u$, $H_2(u) = u^2 - 1$, $H_3(x) = u^3 - 3u$, $H_5(x) = u^5 - 10u^3 + 15u$.

**Remark 2.1.** Let the elements of the population are non-random real numbers, viz. $Y_{mN} = a_{mN}$. Put

$$\bar{a} = N^{-1}(a_{1N} + \ldots + a_{NN}), \quad b^2 = N^{-1}\left((a_{1N} - \bar{a})^2 + \ldots + (a_{NN} - \bar{a})^2\right),$$

$$\tilde{a}_{mN} = (a_{mN} - \bar{a})/b, \quad A_k = N^{-1}\left(\tilde{a}_{1N}^k + \ldots + \tilde{a}_{NN}^k\right). \tag{2.14}$$

Then

$$\mathbb{W}_{3,N}(u) = \Phi(u) - \frac{e^{-u^2/2}}{\sqrt{2\pi}} \left\{ \frac{H_2(u)(1-2p)}{6\sqrt{nq}} A_3 + \frac{H_5(u)(1-4pq)}{72nq} A_3^2 \right.$$

$$\left. + \frac{H_3(u)}{24nq}\left((1-6pq)A_4 - 3(1-4pq)\right) + \frac{H_1(u)p}{2n} \right\}.$$

□

**Remark 2.2**. If $\gamma = 0$ and $\alpha_{20} = 1$ then the ch.f. $\psi_m(t,\tau)$ can be written in the following form

$$\psi_m(t,\tau) = q\exp\left\{-\frac{itpEY_{mN}}{\sqrt{n}\sigma} - \frac{i\tau p}{\sqrt{nq}}\right\} + pE\exp\left\{\frac{it(Y_{mN} - pEY_{mN})}{\sqrt{n}\sigma} + \frac{i\tau q}{\sqrt{nq}}\right\}. \tag{2.15}$$

□

Formula (2.3) with $\psi_m(t,\tau)$ written in the form (2.15) has been used (without proof) by Zhao et al (2004) and Hu et al (2007a). They subsequently proved few lemmas on the integrant of $\Theta_N(t)$ by working with (2.15) in order to get bounds for the remainder terms in the first and the second order asymptotic of the ch.f. $\varphi_n(t)$. In contrast, based on the formula (2.3) and idea that the integrand of $\Theta_N(t)$ is the ch.f. of a sum of $N$ independent two-dimensional r.vec.s.(the fact which was not remarked by the referred authors), the algorithm presented above is more stream and general in constructing asymptotic expansion of any "length". It is also much simpler than that of Mirakhmedov (1983) and Bloznelis (2000).

The formula (2.3) assumes that the r.v.s $\xi_m$ are asymptotically not degenerate, hence the case $n = N$, i.e. the case of sum of independent r.v.s, is excluded. In Zhao et al (2004) and Hu et al (2007a) this limitation was covered by using the approach suggested by von Bahr (1972). Let

$$\phi_m(u) = Ee^{iuY_{mN}}, \quad b_m(t) = e^{t^2\alpha_{20}/2n\sigma^2} Ee^{itY_{mN}/\sigma\sqrt{n}} - 1, \tag{2.16}$$

$$B_j(t) = \frac{(-1)^{j+1}}{j} \sum_{m=1}^{N} b_m^j(t), \quad j = 1,...,n;$$

$$C(n,N,r) = \begin{cases} \dfrac{C_{N-r}^{n-r}}{p^r C_N^n}, & r \leq n \\ 0, & r > n \end{cases}.$$

von Bahr (1972) have proved the following formula

$$e^{t^2\alpha_{20}/2\sigma^2}\varphi_n(t) = \sum_{i_j \geq 0, 1 \leq j \leq n} \prod_{j=1}^{n} \frac{(pB_j(t))^{i_j}}{i_j!} C\left(n, N, \sum_{j=1}^{n} ji_j\right). \tag{2.17}$$

In respect of formula (2.17) we make the following comments. Let us write 1 instead of the $C(n,N,\cdot)$, then the right-hand side has form $\exp\{pB_1(t) + ... + p^n B_n(t)\}$; next using Taylor expansion idea we can observe that the function $p^j B_j(t)$ are polynomial in terms of $\ell_{i,n} t^i / n^{(i-2)/2}$, where $i \geq j$ and $\ell_{i,n}$ are some magnitudes depending on moments of the r.v.s $Y_{m,N}$. Hence, when we wish to approximate $\varphi_n(t)$ by ch.f. of a normal distribution we should approximate $C(n,N,\cdot)$ by 1, separate $\ell_{i,n} t^2$ from $pB_1(t)$ and $p^2 B_2(t)$, and get appropriate bound for $p^3 B_3(t) + ... + p^n B_n(t)$. By this way von Bahr (1972) derived the bound $60 \max_{1 \leq m \leq N} E|Y_{mN} - EY_{mN}|^3 / \sqrt{n}\sigma^3$ for the remainder term in CLT; later Mirakhmedov and Nabiev (1976) had proved the non-uniform variant of von Bahr's result: $c \max_{1 \leq m \leq N} E|Y_{mN} - EY_{mN}|^3 / \sqrt{n}\sigma^3 (1+|x|^3)$; while Zhao et al (2004) obtained Berry-Esseen bound which further improve the result of von Bahr, see below Remark 3.5. When we wish to get $s$-terms Edgeworth asymptotic expansion, in addition to the above procedure, we must separate the terms $\ell_{i,n} t^j / n^{(j-2)/2}$, $j = 2,...,s$ from $pB(t) + ... + p^s B_s(t)$, estimate properly the sum $p^{s+1}B_{s+1}(t) + ... + p^n B_n(t)$, get the asymptotic expansions for $C(n,N,\cdot)$ and $\exp\{pB_1(t) + ... + p^n B_n(t)\}$ using Stirling's formula and Taylor expansion idea respectively. Finally we obtain the asymptotic expansion of $\varphi_n(t)$ by multiplying the asymptotic expansions. This algorithm has been developed by Mirakhmedov (1983) giving the exact formula for the first four terms; the algorithm use more complicated algebra compared to that above presented based on formula (2.3). Recently two terms Edgeworth expansion with improved bound for a remainder term is given by Hu et al (2007a).

In this paper we combine both approaches in developing the methods of aforementioned papers and assuming weaker moment condition.

## 3. Main results.

We still use the notations of Sec1 and 2; see (2.2),(2.5), (2.7), (2.12) and (2.14). Additionally $I\{A\}$ stand for the indicator function of the set $A$ and

$$\Delta_{jN} = \sup_u \left| P\{S_{n,N} < u\sigma\sqrt{n} + n\gamma\} - \mathbb{W}_{j,N}(u) \right|, \ j = 1, 2, 3;$$

where $\mathbb{W}_{1,N}(u) = \Phi(u)$, $\mathbb{W}_{3,N}(u)$ is given in (2.13) and

$$\mathbb{W}_{2,N}(u) = \Phi(u) + \frac{(1-u^2)e^{-u^2/2}}{6\sigma^3\sqrt{2\pi n}}\Lambda_1.$$

**Theorem 3.1**: Let $nq \to \infty$. The Lindberg condition:

$$L_{2N} =: \sum_{m=1}^{N} \tilde{Z}_m^2 I\{|\tilde{Z}_m| \geq \varepsilon\} \to 0 \text{ for arbitrary } \varepsilon > 0, \text{ as } n \to \infty, \tag{3.1}$$

is sufficient for $\Delta_{1N} \to 0$. □

**Remark 3.1.** The condition (3.1) is satisfied if for arbitrary $\varepsilon > 0$

(i) $\quad \dfrac{1}{n\sigma^2} \sum\limits_{m=1}^{N} E(Y_{mN} - pEY_{mN})^2 I\{|Y_{mN} - pEY_{mN}| > \varepsilon\sigma\sqrt{n}\} \to 0,$

and

(ii) $\quad \dfrac{qp}{n\sigma^2} \sum\limits_{m \in D}(EY_{mN})^2 \to 0,$

where $D = \{m : p|EY_{mN}| > \varepsilon\sigma\sqrt{n}\}$, are fulfilled.

Theorem 3.1 gives weakest condition for asymptotic normality of the sample sum from finite population of r.v.s; for the case when $Y_{mN}$ are non-random, the conditions (i) and (ii) coincides with the sufficient and necessary conditions of Erdős Rényi (1959) and Hajek (1960). Also, the condition $nq \to \infty$ is weaker than that assumed by Hajek (1960): $p$ is far away from zero and one. □

Set

$$T_{j,N} = 0.0126 \max\left(\widehat{\beta}_{\min(3,1+j+\delta),N}^{-1}, \left(\beta_{\min(3,1+j+\delta),N} + (nq)^{-1/2}\right)^{-1}\right), \ j = 1, 2, 3, \tag{3.2}$$

where

$$\beta_{kN} = \sum_{m=1}^{N} E|\tilde{Z}_m|^k, \quad \widehat{\beta}_{k,N} = \beta_{k,N}^{(1)} + \sigma^2 \beta_{k,N}^{(2)}, \tag{3.3}$$

$$\beta_{kN}^{(1)} = \frac{N^{-1}\sum\limits_{m=1}^{N} E|Y_{mN} - pEY_{mN} - q\gamma|^k}{n^{(k-2)/2}\sigma^k}, \quad \beta_{kN}^{(2)} = \frac{N^{-1}\sum\limits_{m=1}^{N} |E(Y_{mN} - \gamma)|^k}{n^{(k-2)/2}\sigma^k}, \ k \geq 2.$$

It is easy to see that

$$\beta_{k,N} = \beta_{k,N}^{(1)} + qp^{k-1}\beta_{k,N}^{(2)} \leq \widehat{\beta}_{k,N}. \tag{3.4}$$

**Theorem 3.2**. There exists a constant $C_j > 0$ such that

$$\Delta_{jN} \leq C_j \left( \max\left(\beta_{1+j+\delta,N}, T^{-1}\right) + (nq)^{-(j-1+\delta)/2} + \chi_N\left(T_{j,N}/\sigma\sqrt{n}, T/\sigma\sqrt{n}\right) \right),$$

for any $T \geq T_{j,N}$ and each $j = 1, 2, 3$ where

$$\chi_N(d_0, d_1) = I\{d_0 < d_1\} \int_{d_0 \leq |t| \leq d_1} \left|\frac{\varphi_n(t\sigma\sqrt{n})}{t}\right| dt \leq \min\left(\chi_{1N}(d_0, d_1), \chi_{2N}(d_0, d_1)\right) I\{d_0 < d_1\}, \quad (3.5)$$

$$\chi_{1N}(d_0, d_1) = \sqrt{nq} \ln d_1 \exp\left\{-n\left(1 - \sup_{d_0 \leq |t| \leq d_1} \frac{1}{N} \sum_{m=1}^{N} \left|Ee^{itY_{mN}}\right|\right)\right\},$$

and

$$\chi_{2N}(d_0, d_1) = \sqrt{nq} \ln d_1 \exp\left\{-2nq\left(1 - \sup_{d_0 \leq |t| \leq d_1} \frac{1}{N} \left|\sum_{m=1}^{N} Ee^{itY_{mN}}\right|\right)\right\}. \qquad \square$$

**Remarks 3.2**.

The $\chi_{1N}(d_0, d_1)$ is better then $\chi_{2N}(d_0, d_1)$ under certain conditions such as $q$ is close to zero but many enough of elements of the population are not degenerate r.v.s. In contrast $\chi_{2N}(a,b)$ is applicable for the case when all elements of a population are degenerate, i.e. non-random, viz. $Y_{mN} = a_{mN}$, $m = 1, ..., N$, but $q$ is not too close to zero. In this case $\chi_{2N}(d_0, d_1)$ is exponentially small, for instance, if for a given $d_0 > 0$ there exist $d_1 > 0$, $\varepsilon > 0$ and $\delta > 0$ not depending of $n$ such that

$$\#\left\{k : \left|t(a_{kN} - \bar{a})/b\sqrt{n} - x - 2\pi\nu\right| > \varepsilon\right\} \geq \delta n$$

for any fixed $x$, all $n \geq 1$, $t \in (d_0, d_1)$ and integer $\nu$, see Bickel and van Zwet (1978), Robinson (1978) and Mirakhmedov (1983). $\qquad \square$

Choosing $T = \left(\beta_{2+\delta,N} + (nq)^{-1/2}\right)^{-1}$, $T = \beta_{3+\delta,N}^{-1}$ and $T = \beta_{4+\delta,N}^{-1}$ respectively we obtain from Theorem 2 the following:

**Corollary 3.1**. For any $\delta \in (0,1]$ there exists a positive constants $C_l$ such that

(a) $\quad \Delta_{1N} \leq C_1 \left(\beta_{2+\delta,N} + (nq)^{-\delta/2}\right)$

(b) $\quad \Delta_{2N} \leq C_2 \left(\beta_{3+\delta,N} + (nq)^{-(1+\delta)/2} + \chi_N\left(0.04\sigma^2/V_{3,N}, \left(\beta_{3+\delta,N}\sigma\sqrt{n}\right)^{-1}\right)\right),$

(c) $\quad \Delta_{3N} \leq C_3 \left(\beta_{4+\delta,N} + (nq)^{-(2+\delta)/2} + \chi_N\left(0.04\sigma^2/V_{3,N}, \left(\beta_{4+\delta,N}\sigma\sqrt{n}\right)^{-1}\right)\right),$

where $V_{k,N} = N^{-1} \sum_{m=1}^{N} E\left|Y_{m,N} - \gamma\right|^k$. $\qquad \square$

Taking $T = \widehat{\beta}_{3+\delta,N}^{-1}$ and $T = \widehat{\beta}_{4+\delta,N}^{-1}$, see notation (3.3), we obtain

**Corollary 3.2**. For any $\delta \in (0,1]$ there exists a positive constants $C_l$ such that

(a) $\Delta_{2N} \leq C_2 \left( \hat{\beta}_{3+\delta,N} + (nq)^{-(1+\delta)/2} + \chi_N \left( 0.04\sigma^2 / V_{3,N}, 2^{-(1+\delta)} \sqrt{n^\delta \sigma^{1+\delta}} \right) \right),$

(b) $\Delta_{3N} \leq C_3 \left( \hat{\beta}_{4+\delta,N} + (nq)^{-(2+\delta)/2} + \chi_N \left( 0.04\sigma^2 / V_{3,N}, 2^{-(2+\delta)} \sqrt{n^{1+\delta} \sigma^{2+\delta}} \right) \right).$ □

Set

$$\mu_{k,N} = \frac{p}{(\sigma\sqrt{n})^k} \sum_{m=1}^{N} E|Y_{m,N} - \gamma|^k = V_{k,N} / \sigma^k n^{(k-2)/2}.$$

**Theorem 3.3**. There exists a constant $C_j > 0$ such that for any $\delta \in (0,1]$

$$\Delta_{jN} \leq C_j \left( \max\left(\mu_{1+j+\delta,N}, T^{-1}\right) + \chi_N\left(\tilde{T}_j / \sigma\sqrt{n}, T / \sigma\sqrt{n}\right) \right),$$

for each $j = 1, 2$ and 3, and any $T \geq 0.115 \mu_{\min(3,1+j+\delta),N}^{-1} \stackrel{def}{=} \tilde{T}_j$. □

Choosing $T = 0.115\mu_{2+\delta,N}^{-1}$, $T = \mu_{3+\delta,N}^{-1}$ and $T = \mu_{4+\delta,N}^{-1}$ respectively we obtain from Theorem 3.3 the following

**Corollary 3.3**. For any $\delta \in (0,1]$ there exists a positive constants $C_l$ such that

(a) $\Delta_{1N} \leq C_1 \mu_{2+\delta,N},$

(b) $\Delta_{2N} \leq C_2 \left( \mu_{3+\delta,N} + \chi_N \left( 0.04\sigma^2 / V_{3,N}, 2^{-(1+\delta)} \sqrt{n^\delta \sigma^{1+\delta}} \right) \right),$

(c) $\Delta_{3N} \leq C_3 \left( \mu_{4+\delta,N} + \chi_N \left( 0.04\sigma^2 / V_{3,N}, 2^{-(2+\delta)} \sqrt{n^{1+\delta} \sigma^{2+\delta}} \right) \right).$ □

**Remark 3.3.** We restricted ourselves to just two and three term asymptotic expansions, only to keep the level of complexity of the expressions low. The results can be further extended for $s$-term asymptotic expansion, $s \geq 3$, but at the expense of added complexity in the proof.

**Remarks 3.4**. It is easy to see that $\beta_{k,N} \leq 2^{k-1}(1 + p^{k-1})\mu_{k,N}$. Nevertheless, Theorem 3.3 and Corollary 3.3 gives better bounds for remainder term wrt Theorem 3.2 and Corollary 3.1 in the case when $q$ close to zero; in particular if $q = 0$ then the concerning known results on sum of independent r.v.s. follows from Corollary 3.3.

□

**Remark 3.5.** The main results of Zhao et al (2004) and Hu et al (2007a) respectively state that

$\Delta_{1N} \leq C \min\left( \hat{\beta}_{3,N} + (nq)^{-1/2}, \mu_{3,N} \right),$

$\Delta_{2N} \leq C \min\left( \hat{\beta}_{4,N} + (nq)^{-1}, \mu_{4,N} \right) + \min\left( \chi_{1N} \left( 0.01\sigma^2 / V_{3,N}, 16\sigma\sqrt{n} \right), \chi_{2N} \left( 0.01\sigma^2 / V_{3,N}, 16\sigma\sqrt{n} \right) \right).$

In contrast, Part (a) and (b) of Corollaries 3.1 and 3.3 gives more general bounds using the moments of lower order $2+\delta$ and $3+\delta$, $0 < \delta \leq 1$, respectively. Moreover, since inequality (3.4) the case $\delta = 1$ do improve of their theorems using Lyapunov's ratios $\beta_{k,N}$ instead of $\hat{\beta}_{k,N}$, where $k = 3$ and $k = 4$ respectively. Also Part (c) gives one more term of expansion. If $Y_{mN} = a_{mN}$, $m = 1, ..., N$ are real non-random numbers then the main results of Robinson (1978) and Bloznelis (2000) follows from Corollary 3.2.

□

**Remark 3.6**. It can be readily checked that

$$\left(\alpha_{30} - 3p\alpha_{21} + 2p^2\alpha_{03}\right)/\sigma^3\sqrt{n} = \sum_{m=1}^{N} E\tilde{Z}_m^3. \tag{3.6}$$

Hence the Lyapunov's ratio of the r.v.s $Z_m$, i.e. $\beta_{k,N}$, is natural characteristic in estimate of the remainder terms in approximation by normal distribution, in contrast of $\hat{\beta}_{k,N}$ used by Zhao et al (2004) and Hu et al (2007a). Moreover, due to part (b) of Corollary 3.1 and (3.6) we conjecture that there exist a constant $C > 0$ such that $\Delta_{1N} \leq C\beta_{3,N}$. □

Let $P_N(x) = P\{S_{nN} < x\sigma\sqrt{n} + n\gamma\}$.

**Theorem 3.4.** Let

$$\varlimsup_{n\to\infty} pq > 0 \quad \text{and} \quad \varlimsup_{n\to\infty} N^{-1}\sum_{m=1}^{N} E\exp\{H|Y_{mN}|\} < \infty, \exists\, H > 0. \tag{3.7}$$

Then for all $x \geq 0$, $x = o(\sqrt{N})$

$$\frac{1 - P_N(x)}{1 - \Phi(x)} = \exp\left\{\frac{x^3}{\sqrt{N}} L_N\left(\frac{x}{\sqrt{N}}\right)\right\}\left(1 + O\left(\frac{x+1}{\sqrt{N}}\right)\right), \tag{3.8}$$

$$\frac{P_N(-x)}{\Phi(-x)} = \exp\left\{-\frac{x^3}{\sqrt{N}} L_N\left(-\frac{x}{\sqrt{N}}\right)\right\}\left(1 + O\left(\frac{x+1}{\sqrt{N}}\right)\right), \tag{3.9}$$

where $L_N(v) = \ell_{0N} + \ell_{1N}v + \ldots$ is a power series that for all sufficiently large $N$ is majorized by a power series with coefficients not depending on $N$, and is convergent in some disc, so that $L_N(v)$ converges uniformly in $N$ for sufficiently small values of $|v|$. Particularly

$$\ell_{0N} = \frac{1}{6}\sum_{m=1}^{N} E\tilde{Z}_m^3, \quad \ell_{1N} = \frac{1}{8}\left(\frac{1}{3}\sum_{m=1}^{N} E\tilde{Z}_m^4 - \left(\sum_{m=1}^{N} E\tilde{Z}_m^2\xi_m\right)^2 - \sum_{m=1}^{N}\left(E\tilde{Z}_m^3\right)^2\right). \quad \Box$$

**Corollary 3.4.** Let the conditions (3.7) is satisfied. Then for all $x \geq 0$, $x = o(N^{1/6})$

$$\frac{1 - P_N(x)}{1 - \Phi(x)} = 1 + O\left(\frac{(x+1)^3}{\sqrt{N}}\right) \quad \text{and} \quad \frac{P_N(-x)}{\Phi(-x)} = 1 + O\left(\frac{(x+1)^3}{\sqrt{N}}\right). \tag{3.10}$$

□

**Corollary 3.5.** Let the condition (3.7) is satisfied. Then for all $x \geq 0$, $x = o(N^{1/6})$

$$|P_N(x) - \Phi(x)| = O\left(\frac{(x+1)^2 \exp\{-x^2/2\}}{\sqrt{N}}\right). \tag{3.11}$$

□

**Corollary 3.6.** Let $Y_{mN} = a_{mN}$, $m = 1,\ldots,N$; that is the elements of the population are non-random real numbers. If $\varlimsup_{n\to\infty} pq > 0$ then for all $x \geq 0$, $x = o(\sqrt{N})$ the relations (3.8) and (3.9) hold, and for all $x \geq 0$, $x = o(N^{1/6})$ the relations (3.10) and (3.11) hold. Also, in this case, in notation (2.14)

$$\ell_{0N} = \frac{1-2p}{6\sqrt{pq}} A_3, \quad \ell_{1N} = \frac{1-6pq}{24pq} A_4 - \frac{(1-2p)^2}{8pq} - \frac{1}{8pq} A_3^2.$$

□

**Remark 3.7.** For the case when $Y_{mN} = a_{mN}$, Theorem 1 of Robinson (1977) shows that the relation (3.8) is true for all $x \geq 0$, $x = o\left(\sqrt{N}/\max \tilde{a}_{mN}\right)$. Also we refer to the recently paper of Hu et al (2007b) where it was shown that for

$$0 \leq x \leq C \min\left(\sqrt{nq}/\max|\tilde{a}_{mN}|, (nq)^{1/6}/B_3^{1/3}\right)$$

$$1 - P_N(x) = (1 - \Phi(x))\left(1 + O\left((x+1)^3 B_3/\sqrt{nq}\right)\right),$$

where in addition to the notation (2.14) we denote

$$B_3 = N^{-1}\left(|\tilde{a}_{1N}|^3 + ... + |\tilde{a}_{NN}|^3\right). \tag{3.12}$$

□

**Remark 3.8.** Notice that the sample sum is a special case of the statistics of form $\sum_{m=1}^{N} f_m(\eta_m)$, so namely decomposable statistics (DS), in a simple random sample scheme without replacement, where $f_1, ..., f_N$ are functions defined on the set of non-negative integers, $\eta_m$ is a frequency of *m*-th element of the population in a sample of size $n$. Therefore, from general theorems on DS one can derive some results for the sample sum. In particular, Theorems 3.1 and 3.4 are derived from concerning theorems of Mirakmedov (1996) and Mirakhmedov et al (2012). For details of DS see the referred papers and references within.

□

**4. Lemmas and Proofs**

Before getting to proof of the assertions of Sec 3 we shall prove several auxiliary lemmas. Actually Lemma 4.1 and 4.2 below are valid for any sequence of independent r.vec.s with zero expectations and uncorrelated sums of components; moreover if the r.vec.s has a finite moments of integer order, say $k > 2$, then Lemma 4.1 is a corollary of Theorem 9.11 of BR (1976); authors could not find in literature concerning assertions when the $k$ may be non-integer as we require.

In what follows we still use the notations of Sec 2 and 3, and assume $E|Z_{m,N}|^k < \infty$, where $k > 2$ is not necessary to be integer. Additionally just to keep the notation simple we put $\kappa_{l,N} = \sum_{m=1}^{N} E|\tilde{\xi}_m|^l = (q^{l-1} + p^{l-1})(nq)^{-(l-2)/2}$. We shall use the following relation, see formula (10) of Mirakhmedov (1992): for any complex z, integer $l \geq 0$ and $\delta \in (0,1]$

$$e^z = \sum_{j=0}^{l} \frac{z^j}{j!} + \theta \frac{2^{1-\delta} z^{l+\delta}}{(1+\delta) \cdot ... \cdot (l+\delta)} e^{|\text{Re}\, z|}, \tag{4.1}$$

here and everywhere in what follows $\delta \in (0,1]$, and we use the same symbol $\theta$ for a value such that $|\theta| \leq 1$, although they may not be the same at a different occurrence; also we use well-known inequalities between moments and Lyapunov's ratio: $E|\xi|^j \leq \left(E|\xi|^i\right)^{j/i}$, $1 \leq j \leq i$; $\beta_{l,N}^{1/(l-2)} \leq \beta_{k,N}^{1/(k-2)}$, $3 \leq l \leq k$. We will not refer to these formulas in the cases when it is evident from the context of notations. We assume that, see the denotes (2.2) and (2.12),

$$\gamma = 0 \text{ and } \alpha_{20} = 1, \tag{4.2}$$

the general case can be reduced to that of satisfying (4.2) by considering $(Y_{mN} - \gamma)/\sqrt{\alpha_{20}}$ instead of $Y_{mN}$. So, from now on $Z_m = Y_{mN}\xi_m - pEY_{mN}$, $\sigma^2 = 1 - p\alpha_{02}$, $\mathcal{L}(\xi_m) = Bi(1, p)$. Set

$$\tilde{P}_{1,N}(t,\tau) = 1, \quad \tilde{P}_{2,N}(t,\tau) = 1 + \frac{1}{\sqrt{N}} P_{1,N}(t,\tau), \quad \tilde{P}_{3,N}(t,\tau) = 1 + \frac{1}{\sqrt{N}} P_{1,N}(t,\tau) + \frac{1}{N} P_{2,N}(t,\tau),$$

$$\Im_{j,N} = c_1 \min\left(\beta_{1+j+\delta,N}^{-1/(1+j+\delta)}, \kappa_{1+j+\delta,N}^{-1/(1+j+\delta)}\right).$$

**Lemma 4.1**. If

$$\beta_{1+j+\delta,N} \leq 1 \quad \text{and} \quad \max(|t|,|\tau|) \leq \Im_{j,N}, \tag{4.3}$$

a positive $c_1 \leq 0.1$, then there exist $C_1 > 0$ such that for each $k = 0, 1$ and $j = 1, 2, 3$ one has

$$\left|\frac{\partial^k}{\partial t^k}\left(\prod_{m=1}^{N} \psi_m(t,\tau) - e^{-\frac{t^2+\tau^2}{2}} \tilde{P}_{j,N}(t,\tau)\right)\right|$$

$$\leq C_1 \left(\beta_{1+j+\delta,N} + \frac{1}{(nq)^{(j-1+\delta)/2}}\right)\left(1 + |t|^{1+j+\delta} + |\tau|^{1+j+\delta}\right) e^{-\frac{t^2+\tau^2}{4}}.$$

**Proof**. The case $j = 1$ follows from Part (1) of Lemma A of Mirakhmedov (2005). We restrict ourselves by proof of the case $j = 3$; the case $j = 2$ can be proved in a very similar manner with less algebra. The proof uses an idea of the proof of Theorem 8.6 of BR p.64. We have

$$|\psi_m(t,\tau) - 1| \leq \frac{1}{2} E\left(t\tilde{Z}_m + \tau\tilde{\xi}_m\right)^2 \leq \left(t^2\left(E|\tilde{Z}_m|^{4+\delta}\right)^{2/(4+\delta)} + \tau^2\left(E|\tilde{\xi}_m|^{4+\delta}\right)^{2/(4+\delta)}\right)$$

$$\leq t^2 \beta_{4+\delta,N}^{2/(4+\delta)} + \tau^2 \kappa_{4+\delta,N}^{2/(4+\delta)} < 2c_1^2, \tag{4.4}$$

Hence we can write

$$\prod_{m=1}^{N} \psi_m(t,\tau) = \exp\left\{\sum_{m=1}^{N} \ln \psi_m(t,\tau)\right\} = \exp\left\{-\sum_{m=1}^{N} \sum_{r=1}^{\infty} \frac{1}{r}(1 - \psi_m(t,\tau))^r\right\}.$$

Note that (see (2.8) and (2.9))

$$\frac{1}{\sqrt{N}} \frac{\partial}{\partial t} P_1(t,\tau) = \frac{i^3}{2} \sum_{m=1}^{N} E\tilde{Z}_m\left(t\tilde{Z}_m + \tau\tilde{\xi}_m\right)^2,$$

$$\frac{1}{\sqrt{N}} \frac{\partial}{\partial t} P_1(t,\tau) = \frac{i^4}{6} \sum_{m=1}^{N} E\tilde{Z}_m\left(t\tilde{Z}_m + \tau\tilde{\xi}_m\right)^3 - \frac{i^4}{2} \sum_{m=1}^{N} E\left(t\tilde{Z}_m + \tau\tilde{\xi}_m\right)^2 E\tilde{Z}_m\left(t\tilde{Z}_m + \tau\tilde{\xi}_m\right)$$

$$+ \frac{i^6}{12} \sum_{m=1}^{N} E\left(t\tilde{Z}_m + \tau\tilde{\xi}_m\right)^3 \sum_{m=1}^{N} E\tilde{Z}_m\left(t\tilde{Z}_m + \tau\tilde{\xi}_m\right)^2.$$

To keep the notations simple we set

$$\nabla_N(t,\tau) = \sum_{m=1}^{N}\left((1 - \psi_m(t,\tau)) + 2^{-1}(1 - \psi_m(t,\tau))^2\right) - \frac{t^2 + \tau^2}{2}.$$

Now,

$$\prod_{m=1}^{N} \psi_m(t,\tau) - \exp\left\{-\frac{t^2+\tau^2}{2}\right\}\left(1 + \frac{1}{\sqrt{N}}P_1(t,\tau) + \frac{1}{N}P_2(t,\tau)\right)$$

$$= \exp\left\{-\frac{t^2+\tau^2}{2}\right\}\left[\exp\{-\nabla_N(t,\tau)\} - 1 - \frac{1}{\sqrt{N}}P_1(t,\tau) - \frac{1}{N}P_2(t,\tau)\right]$$

$$+ \exp\left\{-\nabla_N(t,\tau) - \frac{t^2+\tau^2}{2}\right\}\left(\exp\left\{-\sum_{m=1}^{N}\sum_{j=3}^{\infty} j^{-1}(1-\psi_m(t,\tau))^j\right\} - 1\right)$$

$$= \exp\left\{-\frac{t^2+\tau^2}{2}\right\}[A_1(t,\tau) + A_2(t,\tau)]. \tag{4.5}$$

Also

$$\frac{\partial}{\partial t}A_1(t,\tau) = -tA_1(t,\tau) + A_{11}(t,\tau), \tag{4.6}$$

where

$$A_{11}(t,\tau) = \exp\left\{-\frac{t^2+\tau^2}{2}\right\}\left[\exp\{-\nabla_N(t,\tau)\}\left(\sum_{m=1}^{N}(2-\psi_m(t,\tau))\frac{\partial}{\partial t}\psi_m(t,\tau) + t\right.\right.$$

$$\left. -\frac{i^3}{2}\sum_{m=1}^{N}E\tilde{Z}_m\left(t\tilde{Z}_m + \tau\tilde{\xi}_m\right)^2 - \frac{i^4}{6}\sum_{m=1}^{N}E\tilde{Z}_m\left(t\tilde{Z}_m + \tau\tilde{\xi}_m\right)^3 + \frac{i^4}{2}\sum_{m=1}^{N}E\left(t\tilde{Z}_m + \tau\tilde{\xi}_m\right)^2 E\tilde{Z}_m\left(t\tilde{Z}_m + \tau\tilde{\xi}_m\right)\right)$$

$$+ \frac{i^3}{2}\left(\exp\{-\nabla_N(t,\tau)\} - 1 - \frac{i^3}{6}\sum_{m=1}^{N}E\left(t\tilde{Z}_m + \tau\tilde{\xi}_m\right)^3\right)\sum_{m=1}^{N}E\tilde{Z}_m\left(t\tilde{Z}_m + \tau\tilde{\xi}_m\right)^2$$

$$\left. + \frac{i^4}{2}\left(\exp\{-\nabla_N(t,\tau)\} - 1\right)\left(\frac{1}{3}\sum_{m=1}^{N}E\tilde{Z}_m\left(t\tilde{Z}_m + \tau\tilde{\xi}_m\right)^3 - \sum_{m=1}^{N}E\left(t\tilde{Z}_m + \tau\tilde{\xi}_m\right)^2 E\tilde{Z}_m\left(t\tilde{Z}_m + \tau\tilde{\xi}_m\right)\right)\right]$$

$$\equiv A_{11}^{(1)}(t,\tau) + A_{11}^{(2)}(t,\tau) + A_{11}^{(3)}(t,\tau). \tag{4.7}$$

Next,

$$\frac{\partial}{\partial t}A_2(t,\tau) = A_2(t,\tau)\sum_{m=1}^{N}(2-\psi_m(t,\tau))\frac{\partial}{\partial t}\psi_m(t,\tau)$$

$$+ \exp\left\{-\nabla_N(t,\tau) - \frac{t^2+\tau^2}{2}\right\}\exp\left\{-\sum_{m=1}^{N}\sum_{j=3}^{\infty} j^{-1}(1-\psi_m(t,\tau))^j\right\}\sum_{m=1}^{N}\sum_{j=3}^{\infty}(1-\psi_m(t,\tau))^{j-1}\frac{\partial}{\partial t}\psi_m(t,\tau)$$

$$\equiv A_2^{(1)}(t,\tau) + A_2^{(2)}(t,\tau). \tag{4.8}$$

We have

$$\sum_{m=1}^{N}(1-\psi_m(t,\tau)) = \frac{t^2+\tau^2}{2} - \frac{i^3}{6}\sum_{m=1}^{N}E\left(t\tilde{Z}_m + \tau\tilde{\xi}_m\right)^3 + r(t,\tau), \tag{4.9}$$

with $|r(t,\tau)| \leq 3^{-1}(t^4\beta_{4,N} + \tau^4\kappa_{4,N})$, and

$$\sum_{m=1}^{N}|1-\psi_m(t,\tau)|^2 \leq 2(t^4\beta_{4,N} + \tau^4\kappa_{4,N}). \tag{4.10}$$

Using the last and first inequalities of (4.4) we obtain

$$\left|\sum_{m=1}^{N}\sum_{r=3}^{\infty} r^{-1}(1-\psi_m(t,\tau))^r\right| \le \frac{1}{3(1-2c_1^2)} \sum_{m=1}^{N} |1-\psi_m(t,\tau)|^3$$

$$\le \frac{4c_1^{2-\delta}}{3(1-2c_1^2)}\left(|t|^{4+\delta}\beta_{4+\delta,N} + |\tau|^{4+\delta}\kappa_{4+\delta,N}\right) \le \frac{8c_1^6}{3(1-2c_1^2)}. \tag{4.11}$$

Now, use Taylor expansion idea, Eqv. (2.6), the inequalities between Lyapunov's ratios, condition (4.3) and Eqv. (4.9), (4.10) to get (see notation (2.8), (2.9))

$$-\nabla_N(t,\tau) = \frac{i^3}{6}\sum_{m=1}^{N} E\left(t\tilde{Z}_m + \tau\tilde{\xi}_m\right)^3$$

$$+ \frac{i^4}{24}\sum_{m=1}^{N} E\left(t\tilde{Z}_m + \tau\tilde{\xi}_m\right)^4 - \frac{i^4}{8}\sum_{m=1}^{N}\left(E\left(t\tilde{Z}_m + \tau\tilde{\xi}_m\right)^2\right)^2 + r_1(t,\tau)$$

$$\equiv \frac{1}{\sqrt{N}} P_{1,N}(t,\tau) + \frac{1}{N}\left(P_{2,N}(t,\tau) - \frac{1}{2}P_{1,N}^2(t,\tau)\right) + r_1(t,\tau) \tag{4.12}$$

$$= \frac{i^3}{6}\sum_{m=1}^{N} E\left(t\tilde{Z}_m + \tau\tilde{\xi}_m\right)^3 + r_2(t,\tau) \tag{4.13}$$

$$= r_3(t,\tau), \tag{4.14}$$

$$\nabla_N^2(t,\tau) = \frac{1}{N} P_{1,N}^2(t,\tau) + r_4(t,\tau), \tag{4.15}$$

$$\left|\nabla_N(t,\tau)\right|^3 \le r_5(t,\tau), \tag{4.16}$$

where

$$|r_l(t,\tau)| \le C(c_1,\delta)\left(|t|^{4+\delta}\beta_{4+\delta,N} + |\tau|^{4+\delta}\kappa_{4+\delta,N}\right),\ l=1,4,5; \tag{4.17}$$

$$|r_2(t,\tau)| \le 4\left(t^4\beta_{4+\delta,N} + \tau^4\kappa_{4+\delta,N}\right), \tag{4.18}$$

$$|r_3(t,\tau)| \le \frac{2}{3}\left(|t|^3 \beta_{4+\delta,N}^{1/(2+\delta)} + |\tau|^3 \kappa_{4+\delta,N}^{1/(2+\delta)}\right) \le \frac{2c_1}{3}(t^2+\tau^2) < \frac{1}{15}(t^2+\tau^2); \tag{4.19}$$

and

$$\sum_{m=1}^{N}(2-\psi_m(t,\tau))\frac{\partial}{\partial t}\psi_m(t,\tau) = -t + \frac{i^3}{2}\sum_{m=1}^{N} E\tilde{Z}_m\left(t\tilde{Z}_m + \tau\tilde{\xi}_m\right)^2$$

$$+ \frac{i^4}{6}\sum_{m=1}^{N} E\tilde{Z}_m\left(t\tilde{Z}_m + \tau\tilde{\xi}_m\right)^3 - \frac{i^4}{2}\sum_{m=1}^{N} E\left(t\tilde{Z}_m + \tau\tilde{\xi}_m\right)^2 E\tilde{Z}_m\left(t\tilde{Z}_m + \tau\tilde{\xi}_m\right) + r_6(t,\tau), \tag{4.20}$$

where

$$|r_6(t,\tau)| \le 4\left(\left(|t|^{4+\delta} + |\tau|^{4+\delta}\right)\beta_{4+\delta,N} + |\tau|^{4+\delta}\kappa_{4+\delta,N}\right).$$

Also

$$\left|\sum_{m=1}^{N}(2-\psi_m(t,\tau))\frac{\partial}{\partial t}\psi_m(t,\tau)\right|\leq 3(|t|+|\tau|). \tag{4.21}$$

At last, using (4.4) and that

$$\left|\frac{\partial}{\partial t}\psi_m(t,\tau)\right|\leq |t|E\tilde{Z}_m^2+|\tau|E\left|\tilde{Z}_m\xi_m\right|$$

we have

$$\left|\sum_{m=1}^{N}\sum_{j=3}^{\infty}(1-\psi_m(t,\tau))^{j-1}\frac{\partial}{\partial t}\psi_m(t,\tau)\right|\leq \frac{1}{1-2c_1^2}\sum_{m=1}^{N}\left|1-\psi_m(t,\tau)\right|^2\left|\frac{\partial}{\partial t}\psi_m(t,\tau)\right|$$

$$\leq C(c_1,\delta)\left(|t|^{4+\delta}+|\tau|^{4+\delta}\right)\left(\beta_{4+\delta,N}+\kappa_{4+\delta,N}\right). \tag{4.22}$$

To prove Lemma 4.1 it is sufficient to get appropriate upper bound for the module of $A_1(t,\tau)$, $A_2(t,\tau)$ and their derivative. To get a bound

(i) for the $|A_1(t,\tau)|$: use (4.1) with $\delta=1$, $l=2$ and (4.12) - (4.16) and the last inequality of (4.19);

(ii) for the $|A_2(t,\tau)|$: use (4.1) with $\delta=1$, $l=0$ and both inequalities of (4.10), equality (4.13), the last inequality of (4.19);

(iii) for the $|A_{11}(t,\tau)|$: successively to estimate $|A_{11}^{(1)}(t,\tau)|$ use (4.20), (4.14) and (4.19); to estimate $|A_{11}^{(2)}(t,\tau)|$ use (4.1) with $\delta=1$, $l=1$ and (4.15), (4.13), (4.18), (4.14), (4.18); to estimate $|A_{11}^{(3)}(t,\tau)|$ use (4.1) with $\delta=1$, $l=0$ and (4.14) with the first inequality of (4.19);

(iv) for the $\left|\frac{\partial}{\partial t}A_2(t,\tau)\right|$: use to estimate $A_2^{(1)}(t,\tau)$ term the bound already obtained for the $|A_2(t,\tau)|$ and (4.21); next for the $A_2^{(2)}(t,\tau)$ term use (4.22), (4.11), (4.14) and the last inequality of (4.19).

Applying the results of (i)-(iv) in the Eqv.s (4.5), (4.6), (4.7) and (4.8), also using the inequalities between Lypunov's ratios, the condition (4.3) with some algebra one can complete the proof of Lemma 4.1; the details are omitted. □

**Lemma 4.2.** Let $\max\left(\beta_{2+\delta,N},(nq)^{-1/2}\right)\leq 0.01$. If $|t|\leq 0.315\beta_{2+\delta,N}^{-1}$ and $|\tau|\leq 0.315\sqrt{nq}$ then for each $k=0,1$ one has

$$\left|\frac{\partial^k}{\partial t^k}\prod_{m=1}^{N}\psi_m(t,\tau)\right|\leq \exp\left\{-\frac{t^2+\tau^2}{10}\right\}.$$

**Proof of Lemma 4.2** follows from Part (2) of Lemma A of Mirakhmedov (2005). □

**Lemma 4.3**. There exists a $c_1>0$ such that if $|t|\leq c_1\left(\beta_{2+\delta,N}+(nq)^{-1/2}\right)^{-1}$ and $0.0625\sqrt{nq}\leq |\tau|\leq \pi\sqrt{nq}$ then for each $k=0,1$ one has

$$\left|\frac{\partial^k}{\partial t^k}\prod_{m=1}^{N}\psi_m(t,\tau)\right|\leq \exp\{-0.09nq\}.$$

**Proof.** For a given r.v. $\varsigma$ say, let $\varsigma^*=\varsigma-\varsigma'$, where $\varsigma'$ is an independent copy of $\varsigma$. We have

$$|\psi_m(t,\tau)|^2 = |\psi_m(0,\tau)|^2 + E\left[\left(e^{it\tilde{Z}_m^*}-1\right)\left(e^{it\tilde{\xi}_m^*}-1\right)\right] + E\left(e^{it\tilde{Z}_m^*}-1\right)$$

$$\leq |\psi_m(0,\tau)|^2 + |t||\tau|E|\tilde{Z}_m^*\tilde{\xi}_m^*| + t^2 E\tilde{Z}_m^{*2}.$$

Using this and the inequalities $x < \exp\{x-1\}$,

$$\left|\frac{\partial}{\partial t}\prod_{m=1}^{N}\psi_m(t,\tau)\right| \leq \sum_{m=1}^{N}\left|\frac{\partial}{\partial t}\psi_m(t,\tau)\right|\prod_{l\neq m}|\psi_l(t,\tau)|, \quad \sum_{m=1}^{N}\left|\frac{\partial}{\partial t}\psi_m(t,\tau)\right| \leq |t| + 2|\tau|,$$

we find for $k = 0, 1$

$$\left|\frac{\partial^k}{\partial t^k}\prod_{m=1}^{N}\psi_m(t,\tau)\right| \leq \sqrt{e}\left(|t|+|\tau|\right)^k \exp\left\{|t||\tau| + 2t^2 - \frac{1}{2}\sum_{m=1}^{N}\left(1-|\psi_m(0,\tau)|^2\right)\right\}. \tag{4.23}$$

On the other hand $|E\exp\{i\tau\xi_m\}|^2 \leq \exp\{-4pq\sin^2\tau/2\}$,

$$\sin^2\frac{\tau}{2} \geq \frac{\tau^2}{\pi^2}, |\tau| \leq \pi, \quad 1-e^{-u} \geq \frac{1-e^{-c}}{c}u, 0 \leq u \leq c,$$

Hence for $0.0625\sqrt{nq} \leq |\tau| \leq \pi\sqrt{nq}$ we have $\sum_{m=1}^{N}\left(1-|\psi_m(0,\tau)|^2\right) \geq (1-e^{-1})nq/4$. This inequality together with (4.23) yield

$$\left|\frac{\partial^k}{\partial t^k}\prod_{m=1}^{N}\psi_m(t,\tau)\right| \leq \sqrt{e}\left(|t|+|\tau|\right)^k \exp\left\{\pi|t|\sqrt{nq} + 2t^2 - \frac{1-e^{-1}}{4}nq\right\}$$

$$\leq \sqrt{e}\left(|t|+|\tau|\right)^k \exp\left\{-nq\left((1-e^{-1})4^{-1}-(\pi+2c_1)c_1\right)\right\} \leq \exp\{-0.09nq\},$$

since $|t| \leq c_1\left(\beta_{2+\delta,N} + (nq)^{-1/2}\right)^{-1} \leq c_1\sqrt{nq}$, and here we choose $c_1 \leq (1-e^{-1})/16(\pi+1) \approx 0.0126$. Lemma 4.3 follows. □

**Lemma 4.4.** There exists a $c_1 > 0$ such that if $|t| \leq c_1\hat{\beta}_{2+\delta,N}^{-1}$, see (3.3), and $0.0625\sqrt{nq} \leq |\tau| \leq \pi\sqrt{nq}$ then for each $k = 0,1$ one has

$$\left|\frac{\partial^k}{\partial t^k}\prod_{m=1}^{N}\psi_m(t,\tau)\right| \leq 3\left(|t|+|\tau|\right)^k \exp\{-cnq\}.$$

**Proof.** The proof follows by using the Newton-Leibniz rule for derivative of a product, reasons very similar to that of the proof of Lemma 2 of Zhao et al (2004) and inequality

$$\left|\sum_{m=1}^{N}\frac{\partial}{\partial t}\psi_m(t,\tau)\right| \leq |t| + |\tau|. \tag{4.24}$$

Here $c_1 > 0$ so small that $8c_1^2 + c_1 \leq (1-\cos(1/16))/24$, hence $0 < c_1 \leq 0.0329$. □

Let

$$W_{1n}(t) = \exp\{-t^2/2\}, \quad W_{2n}(t) = e^{-t^2/2}\left(1 + \frac{(it)^3}{6\sqrt{n\sigma^3}}\Lambda_1\right), \tag{4.25}$$

and $W_{3n}(t)$ be defined in (2.11). Recall the notation (3.2), (3.3) and note that

$$T_{1,N} = 0.01\max\left(\hat{\beta}_{2+\delta,N}^{-1},\left(\beta_{2+\delta,N}+(nq)^{-1/2}\right)^{-1}\right), \quad T_{2,N}=T_{3,N}=0.01\max\left(\hat{\beta}_{3,N}^{-1},\left(\beta_{3,N}+(nq)^{-1/2}\right)^{-1}\right).$$

**Lemma 4.5**. For each $j=1,2,3$ and $k=0,1$ one has: if $\beta_{1+j+\delta,N}\leq 1$, $\max\left(\beta_{2+\delta,N},(nq)^{-1/2}\right)\leq 0.01$ and $|t|\leq T_{j,N}$ then there exist a constant $C>0$ such that

$$\left|\frac{\partial^k}{\partial t^k}\left(\varphi_N(t)-W_{jn}(t)\right)\right|\leq C\left(\left(1+|t|^{2+j}\right)e^{-t^2/12}\left(\beta_{2+j+\delta,N}+\kappa_{2+j+\delta,N}\right)+e^{-cnq}\right).$$

**Proof**. We assume that $\Im_{j,N}\leq T_{j,N}$, see notation of Lemma 4.1; the case $\Im_{j,N}>T_{j,N}$ is simpler needing considerably less algebra.

Let $|t|\leq \Im_{j,N}$. We have

$$\frac{\partial^k}{\partial t^k}\left(\Theta_N(t)-Q_{iN}(t)\right)=\frac{1}{\sqrt{2\pi}}\int_{|\tau|\leq \Im_{j,N}}\frac{\partial^k}{\partial t^k}\left(\Psi_N(t,\tau)-e^{-\frac{t^2+\tau^2}{2}}\left(1+\tilde{P}_{j,N}(t,\tau)\right)\right)d\tau$$

$$+\frac{1}{\sqrt{2\pi}}\int_{\Im_{j,N}\leq|\tau|\leq 0.315\sqrt{nq}}\frac{\partial^k}{\partial t^k}\Psi_N(t,\tau)d\tau+\frac{1}{\sqrt{2\pi}}\int_{0.315\sqrt{nq}\leq|\tau|\leq \pi\sqrt{nq}}\frac{\partial^k}{\partial t^k}\Psi_N(t,\tau)d\tau$$

$$+\frac{1}{\sqrt{2\pi}}\int_{|\tau|\geq \Im_{j,N}}e^{-\tau^2/2}\frac{\partial^k}{\partial t^k}e^{-t^2/2}\left(1+\tilde{P}_{j,N}(t,\tau)\right)d\tau = D_1+D_2+D_3+D_4. \quad (4.26)$$

Applying Lemma 4.1, 4.2 and 4.3 to estimate $D_1, D_2$ and $D_3$ respectively after simple calculations we obtain

$$\left|D_1+D_2+D_3\right|\leq C\left(\left(1+|t|^{2+i}\right)e^{-t^2/12}\left(\beta_{2+j+\delta,N}+\kappa_{2+j+\delta,N}\right)+e^{-cnq}\right). \quad (4.27)$$

It is readily seen that

$$\left|D_4\right|\leq Ce^{-t^2/2}e^{-\pi^2 nq/4}. \quad (4.28)$$

Now let $\Im_{j,N}\leq|t|\leq T_{j,N}$. Write

$$\frac{\partial^k}{\partial t^k}\left(\Theta_N(t)-Q_{jN}(t)\right)=\frac{1}{\sqrt{2\pi}}\int_{|\tau|\leq 0.315\sqrt{nq}}\frac{\partial^k}{\partial t^k}\Psi_N(t,\tau)d\tau+D_3-\frac{\partial^k}{\partial t^k}Q_{jN}(t)$$

$$=\tilde{D}_2+D_3-\frac{\partial^k}{\partial t^k}Q_{iN}(t). \quad (4.29)$$

Due to definition of $Q_{iN}(t)$ it easy to see that: for $\Im_{j,N}\leq|t|\leq T_{j,N}$

$$\left|\frac{\partial^k}{\partial t^k}Q_{jN}(t)\right|\leq C\left(1+|t|^{2+j}\right)e^{-t^2/4}\left(\beta_{2+j+\delta,N}+\kappa_{2+j+\delta,N}\right). \quad (4.30)$$

Note that $T_{1,N}=0.01\beta_{2+\delta,N}^{-1}$, $T_{2,N}=T_{3,N}=0.01\beta_{3,N}^{-1}$, hence we can use Lemma 4.2 to estimate $\tilde{D}_2$; also to estimate $D_3$ we apply Lemmas 4.3 and 4.4; finally we obtain: if $\Im_{j,N}\leq|t|\leq T_{j,N}$ then

$$\left|\tilde{D}_2 + D_3\right| \leq C\left(\left(\beta_{2+j+\delta,N} + \kappa_{2+j+\delta,N}\right)e^{-t^2/12} + e^{-cnq}\right). \quad (4.31)$$

Thus, from (4.26)-(4.31) we have: if $|t| \leq T_{j,N}$ then for $k=0$ and $k=1$

$$\left|\frac{\partial^k}{\partial t^k}\left(\Theta_N(t) - Q_{jN}(t)\right)\right| \leq C\left(\left(1+|t|^{2+j}\right)e^{-t^2/12}\left(\beta_{2+j+\delta,N} + \kappa_{2+j+\delta,N}\right) + e^{-cnq}\right). \quad (4.32)$$

On the other hand using Stirling's formula it is easy to show that

$$\Theta_N(0) = 1 + \frac{1}{N}G_{2,N}(0) + \frac{2\theta}{(nq)^2}. \quad (4.33)$$

Applying (4.32) and (4.33) in the formula (2.3) and using inequality between Lyapunov's ratios with quite evident calculations we complete the proof of Lemma 4.5. □

**Remark 4.1**. It may remarked that Lemma 4.5 can be further extended for $s$-term asymptotic expansion and any $k=0,1,...,s$, $s \geq 3$, but at the expense of added complexity in the proof. Such an extension of Lemma 4.5 can be used to derive asymptotic expansion of the moments of the sample sum. In particular, one can show that, see notation (2.12),

$$E\left(S_{n,N} - ES_{n,N}\right)^3 = n\Lambda_1 + Cn^{3/2}\sigma^3\left(\beta_{3+\delta,N} + (nq)^{-(1+\delta)/2}\right),$$

$$E\left(S_{n,N} - ES_{n,N}\right)^4 = n^2\sigma^4\left(3 + \frac{\Lambda_2}{n\sigma^4} + 4\frac{pq\alpha_{02}}{n\sigma^2} + C\left(\beta_{4+\delta,N} + (nq)^{-(1+\delta/2)}\right)\right).$$

In the case $Y_{mN} = a_{mN}$ are non-random real numbers these formulas have the following form, see notations (2.14) and (3.12),

$$E\left(S_{n,N} - ES_{n,N}\right)^3 = n(1-2p)A_3 + C\sqrt{nq}\left((1-3pq)A_4 + 1\right),$$

$$E\left(S_{n,N} - ES_{n,N}\right)^4 = (nq)^2\left(3 + \frac{1}{nq}\left((1-6pq)A_4 - 3\right) + \frac{16p}{n} + C\frac{1}{(nq)^{3/2}}\left((p^4+q^4)(1+B_5)\right)\right).$$

An upper bound for $ES_{nN}^4$ has been presented in Lemma 2.6 of Rosen (1967). □

**Lemma 4.6**. Let

$$\mu_{2+\delta,N} \leq \min\left(2^{-1}, 2^{-(2+\delta)/2}(nq)^{-\delta/2}\right). \quad (4.34)$$

Then for $|t| \leq c_1(\delta)\mu_{2+\delta,N}^{-1}$, $c_1(\delta) := \left((1+\delta)(2+\delta)/16\right)^{1/\delta}$, one has

$$\prod_{m=1}^{N}|\psi_m(t,\tau)| \leq \exp\left\{-p^2\frac{t^2}{4}\right\}.$$

**Proof**. Set

$$v_{s,m} = E|Y_{m,N}|^s, \quad \bar{v}_{s,m} = \max(1, v_{s,m}), \quad \vartheta_m^2 = VarY_{mN}, \quad V_{s,N} = N^{-1}\sum_{m=1}^{N}v_{s,m}, \quad \vartheta^2 = N^{-1}\sum_{m=1}^{N}\vartheta_m^2 = 1-\alpha_{02}.$$

Since $\mu_{2+\delta,N} = n^{-\delta/2}\sigma^{-(2+\delta)}V_{2+\delta,N}$ and (4.34) we have $2q^{\delta/(2+\delta)}V_{2+\delta}^{2/(2+\delta)} \leq \sigma^2$. Therefore

$2q\alpha_{02} \leq 2qV_{2+\delta}^{2/(2+\delta)} \leq 2q^{\delta/2+\delta}V_{2+\delta}^{2/(2+\delta)} \leq \sigma^2 = \vartheta^2 + q\alpha_{02}$, hence

$$\sigma^2 \geq \vartheta^2 = \sigma^2 - q\alpha_{02} \geq \sigma^2 - \sigma^2/2 \geq \sigma^2/2. \quad (4.35)$$

Apply (4.1) with $l=2$ to get

$$\left|Ee^{itY_{m,N}/\sigma\sqrt{n}}\right|^2 \le 1 - \frac{t^2\vartheta_m^2}{n\sigma^2} + c(\delta)\frac{|t|^{2+\delta} v_{m,2+\delta}}{(n\sigma)^{2+\delta}}, \quad c(\delta) = 4/(1+\delta)(2+\delta).$$

Using this inequality and that $q^2 + 2pq\left|Ee^{itY_{m,N}/\sigma\sqrt{n}}\right| \le 1 - p^2$ we obtain

$$|\psi_m(t,\tau)|^2 = \left|E\exp\left\{it\frac{Y_{m,N}\xi_m}{\sigma\sqrt{n}} + i\tau\frac{\xi_m}{\sqrt{npq}}\right\}\right|^2 \le \left(q + p\left|Ee^{itY_{m,N}/\sigma\sqrt{n}}\right|\right)^2$$

$$\le 1 - p^2\left(1 - \left|Ee^{itY_{m,N}/\sigma\sqrt{n}}\right|^2\right) \le 1 - p^2\frac{t^2\vartheta_m^2}{n\sigma^2} + c(\delta)p^2\frac{|t|^{2+\delta} v_{m,2+\delta}}{\left(\sigma\sqrt{n}\right)^{2+\delta}}.$$

Hence for $|t| \le c_1(\delta)\mu_{2+\delta,N}^{-1}$, we have

$$\prod_{m=1}^{N} |\psi_m(t,\tau)|^2 \le \exp\left\{-p^2\frac{t^2\vartheta^2}{\sigma^2} + c(\delta)p^2|t|^{2+\delta}\mu_{2+\delta,N}\right\}$$

$$\le \exp\left\{-pt^2\left(\frac{\vartheta^2}{\sigma^2} - c(\delta)c_1(\delta)|t|^\delta \mu_{2+\delta,N}\right)\right\} \le \exp\left\{-p^2t^2\left(\frac{\vartheta^2}{\sigma^2} - \frac{\mu_{2+\delta,N}^{1-\delta}}{4}\right)\right\} \le \exp\left\{-\frac{p^2}{4}t^2\right\},$$

since (4.35) and that $\mu_{2+\delta,N} \le 1$. Lemma 4.6 is proved. $\square$

**Lemma 4.7**. Let $p \ge 1/2$ and $\mu_{1+j+\delta,N} \le 1$, $j=1,2,3$. Then for $|t| \le 0.1\mu_{1+j+\delta,N}^{-1/(1+j+\delta)}$ one has

$$\left|\varphi_n(t) - W_{j,N}(t)\right| \le c_3 |t|^j \left(1 + |t|^{1+\delta}\right)\mu_{1+j+\delta,N}e^{-t^2/4}.$$

where $W_{jn}(t)$ are defined in (2.11) and (4.25).

**Proof.** Proof based on the formula (2.17) of Bahr, hence we use denotes (2.16). Also we shall develop some ideas of Mirakhmedov (1983), Zhao et al (2004) and Hu et al (2007a).

Let $c_0(k,\delta) = 2^{1-\delta}/(1+\delta)\cdot\ldots\cdot(k+\delta)$ and $c_2(k,\delta) = \left(10c_0(k,\delta)(k+1)!\right)^{-1/\delta} \le 0.1$, with $k \ge 2$. Due to (4.2) and inequalities between moments we observe that

$$V_{k+\delta,N} \ge 1. \tag{4.36}$$

Hence

$$\mu_{k+\delta,N}^{-1} = n^{(k-2+\delta)/2}\sigma^{k+\delta}/V_{k+\delta,N} \le n^{(k-2+\delta)/2}\sigma^{k+\delta}. \tag{4.37}$$

If

$$|t| \le c_2(k,\delta)\mu_{k+\delta,N}^{-1/(k+\delta)}, \tag{4.38}$$

then also

$$|t| \le c_2(k,\delta)\sigma n^{(k-2+\delta)/2(k+\delta)} \tag{4.39}$$

and

$$t^2/2n\sigma^2 \le c_2^2(k,\delta)/2n^{2/(k+\delta)} \le 0.005. \tag{4.40}$$

Set $d_{ms} = EY_{m,N}^s$. Applying (4.1) and (4.36)-(4.39) with appropriate $l$ and $k$ respectively we obtain

$$b_m(t) = \frac{itd_{m1}}{\sigma\sqrt{n}} + \frac{t^2}{2n\sigma^2}(\alpha_{20} - d_{m2}) + \frac{(it)^3}{6(\sigma\sqrt{n})^3}(d_{m3} - 3d_{m1}\alpha_{20})$$
$$+ \frac{(it)^4}{24(\sigma\sqrt{n})^4}(d_{m4} - 6d_{m2}\alpha_{20} + 3\alpha_{20}^2) + \theta\frac{t^{4+\delta}\bar{v}_{m,4+\delta}}{6(\sigma\sqrt{n})^{4+\delta}}, \qquad (4.41)$$

$$b_m(t) = \frac{itd_{m1}}{\sigma\sqrt{n}} + \frac{t^2}{2n\sigma^2}(\alpha_{20} - d_{m2}) + \frac{(it)^3}{6(\sigma\sqrt{n})^3}(d_{m3} - 3d_{m1}\alpha_{20}) + \theta\frac{t^{3+\delta}\bar{v}_{m,4+\delta}^{(3+\delta)/(4+\delta)}}{2(\sigma\sqrt{n})^{3+\delta}}, \qquad (4.42)$$

$$b_m(t) = \frac{itd_{m1}}{\sigma\sqrt{n}} + \frac{t^2}{2n\sigma^2}(\alpha_{20} - d_{m2}) + \theta\frac{t^3\bar{v}_{m,2+\delta}}{(\sigma\sqrt{n})^{2+\delta}}, \qquad (4.43)$$

$$b_m(t) = \frac{itd_{m1}}{\sigma\sqrt{n}} + \theta\frac{t^2\bar{v}_{m,2}}{n\sigma^2}, \qquad (4.44)$$

$$b_m(t) = 1.1\theta\frac{t\bar{v}_{m,1}}{\sigma\sqrt{n}}. \qquad (4.45)$$

Let us prove the case $j = 3$. Note that

$$|t|\bar{v}_{m,4+\delta}^{1/(4+\delta)}/\sigma\sqrt{n} \leq 0.2. \qquad (4.46)$$

Using (4.37),(4.38), (4.46), inequality $\bar{v}_{m,l} \leq \bar{v}_{m,k+\delta}^{l/(k+\delta)} \leq \bar{v}_{m,k+\delta}$, $1 \leq l \leq 4+\delta$, and successively (4.41)-(4.45), we obtain

$$pB_1(p) = \frac{(it)^3}{6\sigma^3\sqrt{n}}\alpha_{30} + \frac{(it)^4}{24n\sigma^4}(\alpha_{40} - 3\alpha_{20}^2) + \theta\frac{t^{4+\delta}}{6}\mu_{4+\delta,N}, \qquad (4.47)$$

$$p^2B_2(p) = \frac{t^2}{2\sigma^2}p\alpha_{02} - \frac{(it)^3}{2\sigma^3\sqrt{n}}p\alpha_{21}$$
$$- \frac{(it)^4}{24n\sigma^4}p(4\alpha_{31} - 3\alpha_{20}^2 + 3\alpha_{20}^{(2)} - 12\alpha_{20}\alpha_{02}) + \theta t^{4+\delta}\mu_{4+\delta,N}, \qquad (4.48)$$

$$p^3B_3(p) = \frac{(it)^3}{3\sigma^3\sqrt{n}}p^2\alpha_{03} - \frac{(it)^4}{2n\sigma^4}p^2(\alpha_{20}\alpha_{02} - \alpha_{22}) + \theta t^{4+\delta}\mu_{4+\delta,N}, \qquad (4.49)$$

$$p^4B_4(p) = -\frac{(it)^4}{4n\sigma^4}p^3\alpha_{04} + \theta t^{4+\delta}\mu_{4+\delta,N}, \qquad (4.50)$$

$$|p^lB_l(t)| \leq \left(1.01|t|\mu_{4+\delta,N}^{1/(4+\delta)}\right)^l \leq 1.05|t|^{4+\delta}\mu_{4+\delta,N}(0.101)^{l-5}, \text{ for } 5 \leq l \leq n. \qquad (4.51)$$

For $|t| \leq 0.1\mu_{4+\delta,N}^{-1/(4+\delta)}$ from (4.47)-(4.51) we have

$$\sum_{l=5}^n |p^lB_l(t)| \leq 1.113|t|^{4+\delta}\mu_{4+\delta,N}, \qquad (4.52)$$

and

$$\sum_{l=1}^{4} p^l B_l(p) = \frac{t^2}{2\sigma^2} p\alpha_{02} + \frac{(it)^3}{6\sigma^3\sqrt{n}}\Lambda_1 + \frac{(it)^4}{24n\sigma^4}\Lambda_2 + \theta t^{4+\delta}\mu_{4+\delta,N} \quad (4.53)$$

$$= \frac{t^2}{2\sigma^2} p\alpha_{02} + \theta\frac{t^2}{8}, \quad (4.54)$$

since $\mu_{4+\delta,N} \le 1$ and hence $\mu_{4+\delta,N}^{-1/(4+\delta)} \le \mu_{4+\delta,N}^{-1/(2+\delta)} \le \mu_{4,N}^{-1/2} \le \mu_{3,N}^{-1}$. Thus (4.52) and (4.54) gives

$$\sum_{l=1}^{n}\left|p^l B_l(p)\right| = \frac{t^2}{2\sigma^2} p\alpha_{02} + \theta\frac{t^2}{4}. \quad (4.55)$$

Next we shall use arguments similar as in the relations (3.11)-(3.22) of Mirakhmedov (1983) with application above Eqv.s (4.46)-(4.55). Use the Stirling's formula to get for $0 < r \le n$

$$0 \le C(n,N,r) - 1 + \frac{1-p}{2n}r(r-1) \le 3\frac{1-p^2}{n^2}r^4. \quad (4.56)$$

Now rewrite (2.16) in the form

$$e^{t^2\alpha_{20}/2\sigma^2}\varphi_n(t) = I_1 + I_2 + I_3, \quad (4.57)$$

where

$$I_1 = \sum\prod_{l=1}^{n}\frac{(pB_l(t))^{i_l}}{i_l!}C\left(n,N,\sum_{l=1}^{n}li_l\right),$$

here the summation is over all $i_l \ge 0$, $l=1,2,3,4$ and $i_l > 0$ for at least one $l = 5,\ldots,n$;

$$I_2 = \sum_{\substack{i_l \ge 0, \\ 1 \le l \le 4}}\prod_{l=1}^{4}\frac{(pB_l(t))^{i_l}}{i_l!}\left(C\left(n,N,\sum_{l=1}^{4}li_l\right) - 1 + \frac{1-p}{2n}\sum_{l=1}^{4}li_l\left(\sum_{l=1}^{4}li_l - 1\right)\right),$$

$$I_3 = \sum_{\substack{i_l \ge 0, \\ 1 \le l \le 4}}\prod_{l=1}^{4}\frac{(pB_l(t))^{i_l}}{i_l!}\left(1 - \frac{1-p}{2n}\sum_{l=1}^{4}li_l\left(\sum_{l=1}^{4}li_l - 1\right)\right).$$

Because $C(n,N,r) \le 1$ it follows from (4.52) and (4.55)

$$|I_1| \le \exp\left\{\sum_{l=1}^{4}|p^l B_l(t)|\right\}\left(\exp\left\{\sum_{l=5}^{n}|p^l B_l(t)|\right\} - 1\right) \le c|t|^{4+\delta}\mu_{4+\delta,N}\exp\left\{\frac{t^2}{2\sigma^2}p\alpha_{02} + \theta\frac{t^2}{4}\right\}. \quad (4.58)$$

Using (4.56) we obtain

$$|I_2| \le \frac{1}{n^2}\sum_{\substack{i_l \ge 0, \\ 1 \le l \le 4}}\prod_{l=1}^{4}\frac{|pB_l(t)|^{i_l}}{i_l!}\left(\sum_{l=1}^{4}li_l\right)^4$$

$$\le \frac{c}{n^2}\exp\left\{\sum_{l=1}^{4}|p^l B_l(t)|\right\}\sum_{l=1}^{4}\left(|p^l B_l(t)| + |p^l B_l(t)|^4\right) \le \frac{c}{n^2}t^2(1+t^6)\exp\left\{\frac{t^2}{2\sigma^2}p\alpha_{02} + \theta\frac{t^2}{4}\right\}. \quad (4.59)$$

As in Mirakhmedov (1983, equality (3.17)) we have

$$I_3 = \sum_{\substack{i_l \geq 0, \\ 1 \leq l \leq 4}} \prod_{l=1}^{4} \frac{\left(p^l B_l(t)\right)^{i_l}}{i_l!} - \frac{q}{2n} \left( \sum_{\substack{i_l \geq 0, \\ 1 \leq l \leq 4}} \frac{\left(p^2 B_2(t)\right)^{i_2} (4i_2(i_2-1) + 2i_2)}{i_2!} \prod_{l=1, l \neq 2}^{4} \frac{\left(p^l B_l(t)\right)^{i_l}}{i_l!} \right.$$

$$\left. + \sum_{\substack{i_l \geq 0, \\ 1 \leq l \leq 4}} \prod_{l=1}^{4} \frac{\left(p^l B_l(t)\right)^{i_l}}{i_{ll}!} \left( \sum_{l=1}^{4} li_l \left( \sum_{l=1}^{4} li_l - 1 \right) - (4i_2(i_2-1) + 2i_2) \right) \right)$$

$$= \exp\left\{ \sum_{l=1}^{4} p^l B_l(t) \right\} \left( 1 - \frac{q}{n} p^2 B_2(t) \left( 2p^2 B_2(t) + 1 \right) - \frac{q}{2n} \sum_{l,j=1}^{4} c_{l,j} \left(p^l B_l(t)\right)^{k_l} \left(p^j B_j(t)\right)^{k_j} \right)$$

where $k_l + k_j = 2$, $k_l \geq 0$, $k_2 = 1$ and $c_{i,j}$ are some constants (we do not need to have $c_{i,j}$ in explicit form, although it can be found easily as a long formula). Use (4.47)-(4.50),(4.52), (4.54), inequalities $q\alpha_{02} \leq q\alpha_{20} \leq \sigma^2$, after some algebra we obtain for $|t| \leq 0.1 \mu_{4+\delta,N}^{-1/(4+\delta)}$

$$I_3 = \exp\left\{ \frac{t^2}{2\sigma^2} p\alpha_{02} + \theta \frac{t^2}{8} \right\} \left( 1 + \frac{(it)^3}{6\sigma^3 \sqrt{n}} \Lambda_1 + \frac{(it)^6}{72n\sigma^6} \Lambda_1^2 \right.$$

$$\left. + \frac{(it)^4}{24n\sigma^4} \left( \Lambda_2 + 12qp^2 \alpha_{02}^2 \right) + \frac{pq(it)^2 \alpha_{02}}{2n\sigma^2} + c\theta \left( t^3 + t^{4+\delta} \right) \mu_{4+\delta,N} \right). \tag{4.60}$$

The case $j = 3$ of Lemma 4.7 follows from (4.57) - (4.60).

Proof of the cases $j = 1$ and $j = 2$ is a very similar to proof of the case $j = 3$ with quite evident algebra using Eqv.s (4.37)-(4.40) with $k = 2$, $k = 3$ respectively and inequality $0 \leq C(n, N, r) - 1 \leq n^{-1} r^2$ instead of (4.56). The details are omitted. □

**Proof of Theorem 3.1** follows (see Remark 3.7) in immediate manner from Theorem 4.2 of Mirakhmedov et al (2012) by putting $\omega_m = 1$, $f_{m,N}(x) = Y_{mN} \cdot x$, $m = 1,...,N$, and $\mathcal{L}(\xi_m) = Bi(p)$. □

**Proof of Theorem 3.2.** By Esseen's smoothing lemma and the fact that

$$\left| \varphi_n(t) - W_{jn}(t) \right| \leq |t| \max_{|u| \leq |t|} \left| \frac{\partial}{\partial u} \left( \varphi_n(u) - W_{jn}(u) \right) \right|,$$

we have

$$\Delta_{jn} \leq \frac{1}{\pi} \int_{|t| \leq T} \left| \frac{\varphi_n(t) - W_{jn}(t)}{t} \right| dt + \frac{24}{T\sqrt{2\pi}} \tag{4.61}$$

$$\leq \frac{1}{\pi} \max_{|u| \leq 1} \left| \frac{\partial}{\partial u} \left( \varphi_N(u) - W_{jn}(u) \right) \right| + \frac{1}{\pi} \int_{1 \leq |t| \leq T_{j,N}} \left| \frac{\varphi_n(t) - W_{jn}(t)}{t} \right| dt + \frac{1}{\pi} \int_{T_{j,N} \leq |t| \leq T} \left| \frac{W_{jn}(t)}{t} \right| dt$$

$$+ \frac{1}{\pi} \int_{T_{j,N} \leq |t| \leq T} \left| \frac{\varphi_n(t)}{t} \right| dt + \frac{24}{T\sqrt{2\pi}} = J_1 + J_2 + J_3 + J_4 + \frac{24}{T\sqrt{2\pi}}, \tag{4.62}$$

where $W_{jn}(t)$ are defined in (2.11) and (4.25), and $T_{j,N}$ from Lemma 4.5.

Using Lemma 4.5 and taking into account the exponential factor of $W_{jn}(t)$, we obtain

$$J_1 + J_2 + J_3 \leq C\left(\beta_{2+j+\delta,N} + \frac{1}{(nq)^{(j+\delta)/2}}\right). \tag{4.63}$$

Theorem 3.2 follows from (4.62) and (4.63) because $J_4 = \chi_N\left(T_{j,N}/\sigma\sqrt{n}, T/\sigma\sqrt{n}\right)$. □

In view of equality $Ee^{itY_{m,N}\xi_m + i\tau\xi_m} = 1 + p\left(Ee^{i(\tau + tY_m)} - 1\right)$ we have

$$\left|Ee^{itY_{m,N}\xi_m + i\tau\xi_m}\right| \leq 1 - p\left(1 - \left|Ee^{itY_{mN}}\right|\right) \text{ and } \left|Ee^{itY_{m,N}\xi_m + i\tau\xi_m}\right|^2 = 1 - 2pq\left(1 - E\cos(tY_{m,N} + \tau)\right).$$

Therefore,

$$\prod_{m=1}^{N}\left|Ee^{itY_{m,N}\xi_m + i\tau\xi_m}\right| \leq \exp\left\{-n\left(1 - \frac{1}{N}\sum_{m=1}^{N}\left|Ee^{itY_{m,N}}\right|\right)\right\}, \tag{4.64}$$

and

$$\prod_{m=1}^{N}\left|Ee^{itY_{m,N}\xi_m + i\tau\xi_m}\right| \leq \exp\left\{-pq\sum_{m=1}^{N}\left(1 - E\cos(\tau + tY_{m,N})\right)\right\}$$

$$\leq \exp\left\{-nq\left(1 - \frac{1}{N}\left|\sum_{m=1}^{N}Ee^{i\tau + itY_{m,N}}\right|\right)\right\}. \tag{4.65}$$

On the other hand by formula (2.3)

$$\chi_N\left(T_{j,N}/\sigma\sqrt{n}, T/\sigma\sqrt{n}\right) = \frac{\sqrt{nq}}{\Theta_N(0)}\int_{T_{j,N}/\sigma\sqrt{n}\leq|t|\leq T/\sigma\sqrt{n}}\frac{1}{|t|}\int_{|\tau|\leq\pi}\prod_{m=1}^{N}\left|Ee^{itY_{m,N}\xi_m + i\tau\xi_m}\right|d\tau dt. \tag{4.66}$$

The inequality (3.5) follows from (4.64), (4.65) and (4.66). □

**Proof of Theorem 3.3.** Let $p \leq 1/2$, then $q \geq 1/2$, and taking into account that $\sigma^2 \leq 1$ and (4.28), we have

$$(nq)^{-(k-2)/2} \leq (2/n)^{(k-2)/2} \leq (2/n)^{(k-2)/2}V_{k,N} = (2/n\sigma^2)^{(k-2)/2}\sigma^{(k-2)}V_{k,N} \leq \mu_{k,N}, \text{ therefore}$$

$$\beta_{k,N} + (nq)^{-(k-2)/2} \leq 2^k\mu_{k,N} + (2/n)^{(k-2)/2} \leq (2^k + 1)\mu_{k,N}, \text{ since}$$

$$\beta_{k,N} \leq 2^{k-1}(1 + p^{k-1})\mu_{k,N}. \tag{4.67}$$

Thus if $p \leq 1/2$ or $(nq)^{-(k-2)/2} \leq 2^{k/2}\mu_{k,N}$, then Theorem 3.3 follows from Theorem 3.2.

Let now $p > 1/2$ and $(nq)^{-(k-2)/2} > 2^{k/2}\mu_{k,N}$. In view of (4.61) we have

$$\Delta_{jn} \leq \frac{1}{\pi}\left(\int_{|t|\leq\tilde{T}_2}\left|\frac{\varphi_n(t) - W_{jn}(t)}{t}\right|dt + \int_{\tilde{T}_2\leq|t|\leq T}\left|\frac{W_{jn}(t)}{t}\right|dt + \int_{\tilde{T}_2\leq|t|\leq T_2}\left|\frac{\varphi_n(t)}{t}\right|dt + \int_{T_2\leq|t|\leq T}\left|\frac{\varphi_n(t)}{t}\right|dt\right) + \frac{24}{T\sqrt{2\pi}},$$

where $\tilde{T}_2 = 0.1\mu_{1+i+\delta,N}^{-1/(1+i+\delta)}$ and $T_2 = c_1(\delta)\mu_{2+\delta,N}^{-1}$, $c_1(\delta)$ is a constant defined in Lemma 4.5. Applying here Lemma 4.5, 4.6 and definitions of $W_{j,n}(t)$ after simple calculations we complete the proof of Theorem 3.3. □

**Proof of Corollary 3.1** Note that

$$0.01\left(\beta_{2+\delta,N} + (nq)^{-1/2}\right)^{-1} \leq T_{1,N} \text{ and } T_{2,N} = T_{3,N} \geq 0.01\widehat{\beta}_{3,N}^{-1} \geq 0.001\sigma^3\sqrt{n}/V_{3,N}. \tag{4.68}$$

Therefore Corollary 3.1 follows immediately from Theorem 2.2 by putting $T = 0.01\left(\beta_{2+\delta,N} + (nq)^{-1/2}\right)^{-1}$, $T = \beta_{3+\delta,N}^{-1}$ and $T = \beta_{4+\delta,N}^{-1}$ for the cases (i),(ii) and (iii) respectively. □

**Proof of Corollary 3.2**. Let $\alpha_{02} \leq 1/4$. Then

$$(1-2p\alpha_{02})^2 = (1-p\alpha_{02})^2 - 2p\alpha_{02}(1-2p\alpha_{02}) + p^2\alpha_{02}^2 \geq (1-p\alpha_{02})(1-3p\alpha_{02})$$

$\geq (1-p\alpha_{02})(1-3p/4) \geq \sigma^2/4$. Use this fact and inequality between moments to get, for $k \geq 3$,

$$\widehat{\beta}_{k,N} \geq \beta_{k,N}^{(1)} \geq n^{-(k-2)/2}\sigma^{-k}\left(N^{-1}\sum_{m=1}^{N} E(Y_{m,N} - pEY_{m,N})^2\right)^{k/2}$$

$$\geq n^{-(k-2)/2}\left((1-2p\alpha_{02})/\sigma^2\right)^{k/2} \geq n^{-(k-2)/2}(2\sigma)^{-k/2}.$$

Hence

$$\left(\beta_{k,N}^{(1)}\sigma\sqrt{n}\right)^{-1} \leq 2^{k-2}n^{(k-3)/2}\sigma^{(k-2)/2}, \quad k \geq 3. \tag{4.69}$$

Let $\alpha_{02} > 1/4$. Then $\sigma^2\beta_{k,N}^{(2)} \geq n^{-(k-2)/2}\sigma^{-(k-2)}\alpha_{02}^{k/2}$, hence

$$\left(\sigma^2\beta_{k,N}^{(2)}\right)^{-1}/\sigma\sqrt{n} \leq 2^k n^{(k-3)/2}\sigma^{(k-2)/2}, \quad k \geq 3. \tag{4.70}$$

Put now in Theorem 3.2 $T = \widehat{\beta}_{k,N}^{-1} \leq \left(\beta_{k,N}^{(1)}\right)^{-1} + \left(\sigma^2\beta_{k,N}^{(2)}\right)^{-1}$, use (4.68) and (4.69), (4.70) with $k = 3+\delta$ and $k = 4+\delta$ respectively to complete the proof of Corollary 3.2. □

**Proof of Corollary 3.3**. Part (a) immediately follows from Theorem 3.3 by putting $T = 0.115\mu_{2+\delta,N}^{-1}$. Note that $\max\left(\beta_{k,N}^{(1)}, \sigma^2\beta_{k,N}^{(2)}\right) \leq 2^{k-1}\mu_{k,N}$. Therefore Parts (b) and (c) follows from Theorem 3.3 by putting $T = \mu_{k,N}^{-1}$ and using inequalities (4.68), (4.69) and (4.70) with $k = 3+\delta$ and $k = 4+\delta$ respectively.

□

**Proof of Theorem 3.4** follows (see Remark 3.7) from Theorem 12 of Mirakhmedov (1996) by putting $f_{m,N}(x) = Y_{mN} \cdot x$, $m = 1,...,N$, $\mathcal{L}(\xi_m) = Bi(p)$, and noting that $\sigma^2 \geq q$ is bounded away from zero. □

**Proofs of Corollaries 3.4 and 3.5** follow from Theorem 3.4 immediately, because we can assume, without losing of generality, that $\max|a_{mN}| \leq c$.

□

**Reference**.

Babu, G. J and Singh K. (1985) Edgeworth expansion for sampling without replacement from finite population. J. Multivariate. Anal. 17, 261-278.

Babu, G. J. and Bai, Z. D. (1996). Mixtures of global and local Edgeworth expansions and their applications. J. Multivariate. Anal. 59, 282-307.

von Bahr, B. (1972). On sampling from a finite set of independent random variables. Z.Wahrsch. Verw. Geb. 24, 279–286.

Bhattacharya, R. N. and Ranga Rao, R. (1976). Normal approximation and asymptotic expansions. Wiley, New York.

Bickel, P. J. and van Zwet, W. R. (1978). Asymptotic expansions for the power of distribution-free tests in the two-sample problem. Ann. Statist. 6, 937-1004.

Bikelis, A. (1969). On the estimation of the remainder term in the central limit theorem for samples from finite populations. Studia Sci. Math. Hungar. 4, 345-354 (Russian).

Bloznelis, M. (2000). One and two-term Edgeworth expansion for finite population sample mean. Exact results, I;II. Lith. Math. J. 40, 213-227.; 329-340.

Cochran, W.G. (1963) Sampling Techniques, Wiley, New York.

Erdös, P. and Renyi, A. (1959). On the central limit theorem for samples from a finite population. Fubl. Math. Inst. Hungarian Acad. Sci. 4, 49-61.

Hajek J. (1964). Limiting distributions in simple random sampling from a finite population. Publ.Math.Inst. Hungar Acad. Sci., 5, 361-374

Holst L. (1979). A unified approach to limit theorems for urn models. J. Appl., Probab.,**16**,154-162.

Höglund, T.(1978). Sampling from a finite population. A remainder term estimate. Scand.J. Statistic. 5, 69-71.

Hu, Z., Robinson J. and Wang,Q.(2007a). Edgeworth expansion for a sample sum from a finite set of independent random variables. Electronical J.Probab. 12, 1402-1417.

Hu, Z., Robinson J. and Wang,Q.(2007b). Cramer-type large deviations for samples from a finite population. Ann. Statist. 35, 673-696.

Mirakhmedov, Sh. A. and Nabiev I. (1976). On estimations of the remainder terms in the limit theorems for the sample sums from finite populations of the random variables. Limit theorems and Mathematical Statistics, 1,102-111(Russian).

Mirakhmedov Sh.A. (1979). Asymptotic expansion of the distribution of sample sum from finite population of the independent random variables. Reports of Uzbek Academy of Sciences, 9, 468-471(Russian)

Mirakhmedov, Sh. A. (1983). An asymptotic expansion for a sample sum from a finite sample.Theory Probab. Appl. 28, 492-502.

Mirakhmedov Sh.A. (1985). Estimations of the closeness of the distribution of decomposable statistics in the multinomial scheme. Theory Probabl.& Appl., 30, 175-178

[Mirakhmedov S.M., Jammalamadaka S.R. and Ibragim B.M., (2012) On Edgeworth expansion in generalized urn model. J. Theor. Probabl. DOI 10.1007/s10959-012-0454-z.

Mirakhmedov S.A.,(1996). Limit theorems on decomposable statistics in a generalized allocation schemes. Discrete Math. Appl., 6, 379-404.

Robinson, J. (1977). Large deviation probabilities for samples from a finite population. Ann. Probab. 5, 913–925.

Robinson, J. (1978). An asymptotic expansion for samples from a finite population. Ann. Statist. 6, 1004-1011.

Rosen B. (1967). On the central limit theorem for a class of sampling procedures. Z. Wahrsch. Verw. Geb., 7, 103-115

Wald A. and Wolfowitz J. (1944). Statistical tests based on permutations of the observations. Ann. Math. Stat. 15, 358-372.

Zhao, L.C., Wu, C. Q. and Wang, Q. (2004). Berry-Esseen bound for a sample sum from a finite set of independent random variables. J. Theor. Probabl. 17, 557-572.